\RequirePackage[l2tabu, orthodox]{nag}

\documentclass[12pt]{amsart}
\usepackage{fullpage,url,amssymb,enumerate,colonequals}
\usepackage[all]{xy} 
\usepackage{comment}
\usepackage{graphicx}

\usepackage[OT2,T1]{fontenc}

\usepackage{color}

\def\bA{\mathbb{A}}
\def\cA{\mathcal{A}}

\def\bC{\mathbb{C}}

\def\bE{\mathbb{E}}
\def\cE{\mathcal{E}}

\def\bH{\mathbb{H}}

\def\cM{\mathcal{M}}

\def\cO{\mathcal{O}}
\def\bP{\mathbb{P}}

\def\bQ{\mathbb{Q}}

\def\bR{\mathbb{R}}

\def\bZ{\mathbb{Z}}

\def\barM{\overline{\cM}}

\def\wt{\widetilde}

\DeclareMathOperator{\Aut}{Aut}

\DeclareMathOperator{\ct}{ct}

\DeclareMathOperator{\Gr}{Gr}

\DeclareMathOperator{\Id}{Id}

\DeclareMathOperator{\ord}{ord}

\DeclareMathOperator{\pt}{pt}

\DeclareMathOperator{\SL}{SL}
\DeclareMathOperator{\Mp}{Mp}
\DeclareMathOperator{\Sp}{Sp}

\DeclareMathOperator{\Stab}{Stab}
\DeclareMathOperator{\Sym}{Sym}

\DeclareMathOperator{\vir}{vir}

\newcommand{\A}{\mathcal A}
\newcommand{\CH}{\mathrm{CH}}
\newcommand{\Mod}{\mathrm{Mod}}

\newcommand{\bs}{\backslash}

\newcommand{\NL}{{\mathsf{NL}}}
\newcommand{\wNL}{{\widetilde{\mathsf{NL}}}}
\newcommand{\Tor}{{\mathsf{Tor}}}
\newcommand{\km}{\mathrm{km}}

\newtheorem{thm}{Theorem}[section]
\newtheorem{lem}[thm]{Lemma}
\newtheorem{cor}[thm]{Corollary}
\newtheorem{prop}[thm]{Proposition}

\newtheorem{conj}{Conjecture}
\newtheorem{question}{Question}

\theoremstyle{definition}
\newtheorem{rem}[thm]{Remark}

\newtheorem{defn}[thm]{Definition}

 \makeatletter
\let\@wraptoccontribs\wraptoccontribs
\makeatother

\makeatletter
\g@addto@macro\bfseries{\boldmath} 
\makeatother

\usepackage{microtype}

\newcommand{\Z}{\mathbb{Z}}
\newcommand{\C}{\mathbb{C}}
\newcommand{\Q}{\mathbb{Q}}
\DeclareMathOperator{\Ind}{Ind}
\newcommand{\R}{\mathbb{R}}
\DeclareMathOperator{\GL}{GL}
\newcommand{\set}[1]{\left\{#1\right\}}

\DeclareMathOperator{\SO}{SO}
\DeclareMathOperator{\GSp}{GSp}

\renewcommand{\d}{\mathrm{d}}


\begin{document}

\title{Modularity of $d$-elliptic loci with level structure}

\author{Fran\c cois Greer}
\address{Michigan State University, Department of Mathematics, 619 Red Cedar Rd
\hfill \newline\texttt{}
 \indent East Lansing, MI 48824}
 \email{{\tt greerfra@msu.edu}}

\author{Carl Lian}

\address{Tufts University, Department of Mathematics, 177 College Ave
\hfill \newline\texttt{}
 \indent Medford, MA 02155} 
 \email{{\tt Carl.Lian@tufts.edu}}
 \contrib[with an appendix by]{Naomi Sweeting}
\address{Princeton University, Department of Mathematics, Fine Hall, Washington Road \hfill \newline\texttt{}
 \indent Princeton, NJ 08540} 
 \email{{\tt naomiss@princeton.edu}}
\date{\today}

\begin{abstract}
We consider the generating series of special cycles on $\A_1(N)\times \A_g(N)$, with full level $N$ structure, valued in the cohomology of degree $2g$. The modularity theorem of Kudla-Millson for locally symmetric spaces implies that these series are modular. When $N=1$, the images of these loci in $\cA_g$ are the $d$-elliptic Noether-Lefschetz loci, which are conjectured to be modular. In the appendix, it is shown that the resulting modular forms are nonzero for $g=2$ when $N\geq 11$ and $N\neq 12$.
\end{abstract}

\maketitle

\section{Introduction}

\subsection{$d$-elliptic loci}

For integers $d,g\ge1$, let $\wNL_{g,d}$ be the moduli space of morphisms of abelian varieties $h:E\to A$, where $E$ is an elliptic curve, $A$ is a principally polarized abelian variety (PPAV) of dimension $g$ with polarization $\Theta_A$, and $\deg(h^{*}\Theta_A)=d$. Let $[\wNL_{g,d}]\in \CH^{g-1}(\cA_g)$ be the Noether-Lefschetz cycle class associated to the forgetful morphism $\epsilon:\wNL_{g,d}\to \cA_g$. (We work throughout with Chow and cohomology groups with rational coefficients.) Set also $[\wNL_{g,0}]=\frac{1}{24}(-1)^g\lambda_{g-1}\in \CH^{g-1}(\cA_g)$, see \S\ref{d=0_sec}. These classes have attracted much recent attention; see \cite{pand_notes_24} for a survey. The purpose of this paper is to advance the following conjecture.

\begin{conj}\label{conj_ag}
    The generating series
    \begin{equation*}
        \sum_{d\ge 0}[\wNL_{g,d}]q^d \in  \CH^{g-1}(\cA_g)\otimes \bQ[[q]]
    \end{equation*}
    is a cycle-valued modular form of weight $2g$.
\end{conj}

Iribar L\'{o}pez \cite[Corollary 4]{il} shows that Conjecture \ref{conj_ag} is true upon projection to the tautological ring $R^{g-1}(\cA_g)$, in the sense of \cite{cmop}. Further plausibility checks for Conjecture \ref{conj_ag} are afforded by pulling back by the Torelli map $\Tor:\cM_g^{\ct}\to\cA_g$. It is proven in \cite{gl} that the pullbacks $\Tor^{!}[\wNL_{g,d}]$ coincide with the Gromov-Witten virtual classes of the loci of curves $[C]\in \cM_g^{\ct}$ admitting a stable map $f:C\to E$ of degree $d$ to \emph{some} genus 1 curve $E$. In particular, Conjecture \ref{conj_ag} would imply that the generating series
    \begin{equation*}
        \sum_{d\ge 0}[\cM^{\ct,q}_{g,1}(\cE,d)]^{\vir} q^d \in  \CH^{g-1}(\cM_g^{\ct})\otimes \bQ[[q]]
    \end{equation*}
    is a cycle-valued modular form of weight $2g$, where $\cM^{\ct,q}_{g,1}(\cE,d)$ is the global moduli space of pointed stable maps $f:(C,q)\to (E,p)$ from a compact type curve of genus $g$ to a varying elliptic curve. In some cases, the Torelli pullback $\Tor^{!}[\wNL_{g,d}]$ may be understood more explicitly, see \cite{cop}.

The Gromov-Witten classes of stable curves admitting a cover of a \emph{fixed} elliptic curve are shown to be {\it quasi}-modular in \cite{obpix}. On the other hand, the closely related (but different) loci of curves $[C]\in \barM_g$ admitting an \emph{admissible cover} of some elliptic curve $f:C\to E$ of degree $d$ are conjectured to be quasi-modular in \cite{lian_d-ell}, and shown to be so when $g\le 3$. It is furthermore shown in \cite{lian_nontaut} that a certain obstruction to the admissible cover loci being tautological is modular in $d$, for any $g$.

We prove the following result in this paper.

\begin{thm}\label{thm_trivial}
        The generating series
    \begin{equation*}
        \sum_{d\ge 0}[\wNL_{g,d}]^+q^d \in  H^{2g}(\cA_1\times \cA_g)\otimes \bQ[[q]]
    \end{equation*}
    is a cycle-valued modular form of weight $2g$. 
\end{thm}
Here, we consider cycle classes of the maps $\wNL_{g,d}\to \cA_1\times \cA_g$, remembering both $E$ and $A$, rather than only $\wNL_{g,d}\to \cA_g$. We use the notation $[\wNL_{g,d}]^+$ to distinguish from the classes $[\wNL_{g,d}]$ on $\cA_g$. We also set $[\wNL_{g,0}]^+=0$. Theorem \ref{thm_trivial} does not imply Conjecture \ref{conj_ag} (even in cohomology), because there is no proper pushforward map $H^{2g}(\cA_1\times \cA_g)\to H^{2(g-1)}(\cA_g)$. In fact, as stated, Theorem \ref{thm_trivial} is trivial: as we prove in Proposition \ref{N=1_trivial}, the classes $[\wNL_{g,d}]^+$ are zero in $H^{2g}(\cA_1\times \cA_g)$, and are moreover zero in $CH^g(\cA_g\times \cA_1)$ (Remark \ref{chow_vanish_remark})!

To obtain a non-trivial statement, we add level structure to the moduli problem. Our main result is the following.
\begin{thm}\label{thm_levelN}
Fix an integer $N\ge1$ and a symplectic group homomorphism $b:(\bZ/N\bZ)^2\to (\bZ/N\bZ)^{2g}$. Let $\wNL^b_{g,d}(N)$ be the moduli space of morphisms $h:E\to A$ as before, where in addition $E,A$ are endowed with full level-$N$ structure and the map induced on $N$-torsion by $h$ is given by $b$.

Then, the generating series
    \begin{equation*}
        \sum_{d\ge 0}[\wNL^b_{g,d}(N)]^+q^d \in  H^{2g}(\cA_1(N)\times \cA_g(N))\otimes \bQ[[q]]
    \end{equation*}
    is a cycle-valued modular form of weight $2g$ and level $N$.
\end{thm}

We again set $[\wNL^b_{g,0}]^+=0$. In contrast to the situation of Theorem \ref{thm_trivial}, the classes $[\wNL^b_{g,d}(N)]^+$ are proven not all to vanish when $g=2$ and $N=11$ and $N\ge13$ in the appendix.

We will see in Proposition \ref{N=1_trivial} that the classes $[\wNL^b_{g,d}(N)]^+$ are supported in the odd K\"{u}nneth component $H^1(\cA_1(N))\otimes H^{2g-1}(\cA_g(N))$. Thus, Theorem \ref{thm_levelN} witnesses modularity of Noether-Lefschetz cycles in a different part of cohomology as Iribar L\'{o}pez's result, which lives in the tautological (hence even) part.

Theorem \ref{thm_levelN} is proven by expressing the cycles $[\wNL^b_{g,d}(N)]^+$ as pullbacks of special cycles from certain (non-algebraic) symmetric spaces, which we discuss in the next section.

\subsection{Symmetric spaces}
Let $(\Lambda,\omega)$, resp. $(\Lambda',\omega')$ be $\bZ^2$, resp. $\bZ^{2g}$, equipped with the standard symplectic forms. The tensor product $L=\Lambda\otimes \Lambda'$ has a natural symmetric bilinear pairing $\gamma$ given by
\begin{equation*}
    \gamma(v_1\otimes v'_1, v_2\otimes v'_2):= \omega(v_1,v_2)\,\omega'(v'_1,v'_2). 
\end{equation*}
As an integral lattice, we have $L\simeq U^{\oplus 2g}$, where $U$ is the hyperbolic plane lattice with Gram matrix
\begin{equation*}
\begin{pmatrix}  0 & 1 \\ 1 & 0 \end{pmatrix}.
\end{equation*}

Let $V=L\otimes \bR$, and let $\Gamma_L\subset SO(V)$ be the subgroup of integral isometries of $L$ that act trivially on $L^\vee/V$, where $L^\vee$ is the dual lattice to $L$. If $K_\infty\subset SO(V)$ is a maximal compact subgroup, then the double quotient
\begin{equation*}
\km(L) := \Gamma_L\bs SO(V)/K_\infty
\end{equation*}
is an example of a non-compact locally symmetric space studied by Kudla-Millson in \cite{km}. In particular, $\km(L)$ contains a countable collection of totally geodesic cycles $C_d$ indexed by positive integers $d$. Since the lattice $L$ has signature $(2g,2g)$, these are cycles of real codimension $2g$. The main theorem of \cite{km} implies that their Poincar\'{e} duals $[C_d]\in H^{2g}(\km(L),\bQ)$ are the Fourier coefficients of a classical modular form of weight $2g$, level 1:
\begin{thm}[Kudla-Millson \cite{km}]
The generating series
\begin{equation*}
  \Phi(q)=\sum_{d\geq 0} [C_d]q^d \in \Mod(2g, \SL_2(\bZ))\otimes H^{2g}(\km(L)).
\end{equation*}
\end{thm}
By convention, we set $[C_0]:=e(\Lambda_g^\vee)\in H^{2g}(\km(L))$, where $\Lambda_g$ is the tautological real vector bundle of rank $2g$ on $\km(L)$ and $e(\Lambda_g^\vee)$ is the Euler class of its dual. Recall that we work throughout with rational coefficients.

For $g>1$, $\km(L)$ is not an algebraic variety, but by the tensor product construction above, it receives a map
\begin{equation*}
\phi: \A_1\times \A_g \to \km(L),
\end{equation*}
defined in \S\ref{km_sec}. We show in Proposition \ref{special_cycle_pullback} that the classes $[C_d]\in H^{2g}(\km(L),\bQ)$ pull back under $\phi$ to the Noether-Lefschetz cycles $[\wNL_{g,d}]^+\in H^{2g}(\cA_1\times \cA_g)$, giving Theorem \ref{thm_trivial}.

However, as we have already mentioned, we prove in Proposition \ref{N=1_trivial} that in fact the classes appearing in Theorem \ref{thm_trivial} are all zero. To obtain a non-trivial $q$-series, we need to add level structure to the moduli spaces involved. 

Let $N\ge 1$ be an integer, let $\cA_1(N)$ and $\cA_g(N)$ are the moduli spaces of elliptic curves and PPAVs with full level-$N$ structure. As in Theorem \ref{thm_levelN}, let $\wNL^b_{g,d}(N)$ be the moduli space of maps $f:E\to A$ of degree $d$ whose induced map on $N$-torsion is given by a specified matrix $b$. Let $\km(L(N))$ be the Kudla-Millson space defined by replacing $\Gamma_L$ in the definition of $\km(L(N))$ with the subgroup of isometries that reduce to the identity mod $N$. We have an enhanced map
\begin{equation*}
\phi(N): \A_1(N)\times \A_g(N) \to \km(L(N)).
\end{equation*}
Then, by the theorem of Kudla-Millson \cite{km}, we have a generating series
\begin{equation*}
    \Phi^{b}(q):=\sum_{d\ge0}[C^b_d(N)]q^d \in \Mod(2g,\Gamma(N))\otimes H^{2g}(\km(L(N))
\end{equation*}
lifting $\Phi(q)$. The $d$-th Fourier coefficient of $\Phi(q)$ pulls back to the Noether-Lefschetz cycle $[\wNL^b_{g,d}(N)]^+\in H^{2g}(\cA_1(N)\times\cA_g(N))$. The difference here is that the modular curve $\A_1(N)$ has non-trivial first cohomology group. We restate Theorem \ref{thm_levelN} as follows.
\begin{thm}
    The pullback 
    \begin{equation*}
        \phi(N)^* \Phi^{b}(q)=\sum_{d\ge 0}[\wNL^b_{g,d}(N)]^+q^d
    \end{equation*}
    is a modular form of weight $2g$, level $N$, valued in $H^{2g}(\A_1(N)\times \A_g(N))$ with support in the odd Kunneth component
    $$H^1(\A_1(N))\otimes H^{2g-1}(\A_g(N)).$$
\end{thm}

\subsection{Further directions}

In order to gain access to a proper pushforward map relating classes on $\cA_1\times\cA_g$ (possibly with level structure added) to those on $\cA_g$, one needs to add cusps to $\cA_1$. We take here $N=1$ for ease of notation, so that $\cA_1^{*}=\barM_{1,1}$. Then, the natural map $\wNL_{g,d}\to \cA_1^{*}\times\cA_g$ remains proper, because a PPAV contains no rational curves. Thus, we may consider the class $[\wNL_{g,d}]^+\in \CH^{g}(\cA_1^{*}\times\cA_g)$, which pushes forward to the class $[\wNL_{g,d}]\in \CH^{g-1}(\cA_g)$ appearing in Conjecture \ref{conj_ag}. We refine the conjecture as follows:

\begin{conj}\label{quasimod}
    The classes of the compactified $d$-elliptic cycles $[\wNL_{g,d}]^+\in CH^{g}(\A_1^* \times \A_g)$ are the Fourier coefficients of a modular form of weight $2g$, level 1.
\end{conj}

After pushforward, Iribar L\'{o}pez's tautological projection calculation \cite{il} shows that the classes $[\wNL_{g,d}]\in\CH^{g-1}(\cA_g)$ are non-zero in general, in contrast to the situation on $\cA_1\times\cA_g$. It follows that the classes $[\wNL_{g,d}]^{+}\in \CH^{g}(\cA_1^{*}\times\cA_g)$ are non-zero, and because the tautological subspace of $\CH^{g-1}(\cA_g)$ maps injectively to $H^{2(g-1)}(\cA_g)$, also that the classes $[\wNL_{g,d}]^{+}\in H^{2g}(\cA_1^{*}\times\cA_g)$ are non-zero. 

It is natural to expect a passage from modularity to \emph{quasi-}modularity upon extending the cycles $\wNL_{g,d}$ to compactifications of $\cA_g$, consistent with results \cite{egt,gar} establishing such phenomena at the boundary of orthogonal type Shimura varieties. This is also consistent with the calculations in \cite{lian_d-ell,obpix} finding cycled-valued quasi-modular forms on $\barM_g$.

Note that Conjectures \ref{conj_ag} and \ref{quasimod} are formulated with values in the Chow group, following \cite{pand_notes_24,il}. The methods we employ here are more likely to prove the cohomological version, since the space $\km(L)$ is non-algebraic. One can also formulate both of these conjectures with level structure in the obvious way.

Extending the results of \cite{gl} to take level structure into account, the classes $[\NL^b_{g,d}(N)]^+\in H^{2g}(\cA_1(N)\times\cA_g(N))$ pull back under the pointed Torelli map $\Tor_1(N):\cM_{g,1}^{\ct}(N)\to \cA_g(N)$ to the virtual (in the sense of Gromov-Witten theory) loci of curves $C$ admitting a cover $f:C\to E$ inducing the map $b$ on $N$-torsion. After capping with the $\psi$ class on $\cM_{g,1}^{\ct}(N)$ and pushing forward to $\cM^{\ct}_g(N)$, we obtain a modular series of classes in $H^{2g}(\A_1(N)\times \cM_g^{ct}(N),\bQ)$, again supported on the odd Kunneth component
\begin{equation*}
\Mod(2g,\Gamma(N))\otimes H^1(\A_1(N))\otimes H^{2g-1}(\cM_g^{ct}(N)).
\end{equation*}
Dually, we obtain a map $\Mod(2g,\Gamma(N))^* \otimes H_1(\A_1(N)) \to H^{2g-1}(\cM^{\ct}_g(N))$. 

\begin{question}
    What is the image of this map?
\end{question}

\subsection{Acknowledgments} 

We thank Rahul Pandharipande, Samir Canning, and Aitor Iribar L\'{o}pez for helpful discussions related to this paper, and for sharing their conjectures. We also thank the anonymous referee for their helpful comments. F.G. was supported by NSF grant DMS-2302548. C.L. and N.S. were supported by NSF postdoctoral fellowship grants DMS-2001976 and DMS-2401823, respectively.

\section{Lattices and theta functions}
We follow the exposition in \cite{mp} for the definitions of metaplectic groups and vector-valued modular forms.
Let $(L,(,))$ be a positive definite, even integral lattice of rank $r$. The Poisson summation formula implies that 
$$\Theta_L(q) := \sum_{x\in L^\vee} q^{\frac{1}{2}(x,x)}e_{[x]}\in \Mod(r/2,\Mp_2(\bZ),\bQ[L^\vee/L]),$$
where $\Mp_2(\bZ)$ is the integer metaplectic group, acting on $\bQ[L^\vee/L]$ via the Weil representation. This representation factors through a double cover of $\SL_2(\bZ/N\bZ)$, where $N$ is the smallest positive integer such that $N\cdot (x,x)\in 2\bZ$ for all $x\in L^\vee$. In particular, for any fixed coset $\delta\in L^\vee/L$, we have:
$$\Theta_{L,\delta}(q) := \sum_{x\in \delta} q^{\frac{1}{2}(x,x)} \in \Mod(r/2,\Gamma(N)).$$
When $r\in 2\bZ$, this gives a classical modular form of weight $r/2$, level $\Gamma(N)$.

The cohomological theta correspondence of Kudla-Millson allows us to reformulate this story in the setting of the locally symmetric space associated to an indefinite lattice. Let $(L,\gamma)$ be an indefinite, even integral lattice of signature $(r_+,r_-)$ and rank $r$, and set $V = L\otimes \bR$. Let $\Gamma_L\subset SO(V)$ be the subgroup of integral isometries acting trivially on $L^\vee/L$, and let $K_\infty\subset SO(V)$ be a maximal compact subgroup, isomorphic to $SO(r_+)\times SO(r_-)$.
\begin{defn}
    Let $\km(L)= \Gamma_L \bs SO(V)/K_\infty$. Note that the symmetric space $SO(V)/K_\infty$ is naturally identified with $\Gr^-(r_-,V)$, the space of oriented negative definite real $r_-$-planes in $V$.
\end{defn}
In the symmetric space $SO(V)/K_\infty$, we have an infinite arrangement of totally geodesic submanifolds indexed by dual lattice vectors $v\in V^\vee$ with $\gamma(v,v)>0$.
\begin{defn}
    Let $C_v\subset \Gr^-(r_-,V)$ be the set of negative definite $r_-$-planes that are orthogonal to $v$. It is a symmetric subspace isomorphic to $\Gr^-(r_-,v^\perp\otimes\bR)$.
\end{defn}
For each integer $d\geq 0$, and $\delta\in V^\vee/V$, the arithmetic subgroup $\Gamma_L$ acts on the set $\{v\in V^\vee: \gamma(v,v)=2d, [v] = \delta\}$ with finitely many orbits; see Lemma \ref{borel}. This allows one to define finite type cycles in the quotient $\km(L)$.
\begin{defn}
    For each $d>0$ and $\delta\in V^\vee/V$, let $C_{d,\delta}\subset \km(L)$ be the image of 
    $$\bigcup_{\substack{\gamma(v,v)=d \\ [v] = \delta}} C_v \subset \Gr^-(r_-,V)$$
    under the quotient map to $\Gamma_L \bs \Gr^-(r_-,V)$. We define smooth uniformizations of $C_{d,\delta}$ using a finite set of $\Gamma_L$-orbit representatives $v_1,\dots,v_{m(d)}$ among the dual lattice vectors of norm $d$ and class $\delta$:
    $$\widetilde{C}_d := \bigsqcup_{i=1}^{m(d)} C_{v_i}$$
\end{defn}

\begin{thm} \cite{km}
    For each $\delta\in L^\vee/L$, the power series
    $$\Phi_{L,\delta}(q) := e_0+\sum_{d>0} [C_{d,\delta}] q^d \in \Mod(r/2, \Gamma(N))\otimes H^{r_-}(\km(L),\bQ)$$
is a modular form, where $e_0$ is the Euler class of the dual tautological bundle of $r_-$-planes.
\end{thm}
In this paper, we specialize to the case where $L = U^{\oplus 2g}$, where $U$ is the hyperbolic plane lattice. This lattice is unimodular, so $L^\vee/L = \{0\}$. More generally, we will consider $L(N)$, for some level $N>0$, which is the lattice $L$ with the quadratic form values multiplied by $N$.
\begin{prop}
    $$L(N)^\vee/ L(N) \simeq L/NL\simeq (\bZ/N\bZ)^{4g}.$$
\end{prop}
\begin{proof}
Multiplication by $1/N$ induces horizontal isomorphisms of the abelian groups in the following commutative diagram
$$\xymatrix{
L \ar[r]^{\sim\,\,\,} & L(N)^\vee \\
NL \ar[r]^{\sim\,\,\,} \ar[u] & L(N).\ar[u]
}$$
The vertical arrows are inclusions of lattices in the same quadratic space.
Since $L$ has rank $4g$, the discriminant formula follows.
\end{proof}

\section{Noether-Lefschetz loci}

Fix again an integer $N\ge1$.

\begin{defn}
    Let $\cA_g(N)$ be the moduli space of triples $(A,\Theta_A,\iota_A)$, where $(A,\Theta_A)$ is a principally polarized abelian variety (PPAV) of dimension $g$ and $\iota_A:A[N]\to (\bZ/N\bZ)^{2g}$ is a symplectic isomorphism.

\end{defn}

When $g=1$, an elliptic curve is canonically polarized, so we drop $\Theta$ from the notation.


\begin{defn}
    Let $d,g\ge1$ be integers. Let $\wNL^b_{g,d}(N)$ be the moduli space of isomorphism classes of maps of abelian varieties $h:E\to A$, where:
    \begin{itemize}
        \item $(E,\iota_E)\in \cA_{1}(N)$ and $(A,\Theta_A,\iota_A)\in \cA_{g}(N)$,
        \item $h^{*}\Theta_A$ has degree $d$ on $E$, 
        \item $b:(\bZ/N\bZ)^{2} \to (\bZ/N\bZ)^{2g}$ is a fixed group homomorphism respecting the standard symplectic forms, and
        \item the diagram
        \begin{equation*}
        \xymatrix{
        E[N] \ar[r]^{h} \ar[d]^{\iota_E} & A[N] \ar[d]^{\iota_A}\\
        (\bZ/N\bZ)^{2} \ar[r]^b & (\bZ/N\bZ)^{2g}
        }
        \end{equation*}
        commutes.
    \end{itemize}
    An isomorphism of maps $h:E\to A$ and $h':E'\to A'$ consists of the data of isomorphisms $f_E:E\to E'$ and $f_A:A\to A'$, such that $\iota_{E'}\circ f_E=\iota_E$, $\iota_{A'}\circ f_A =\iota_A$, and $f_A^{*}\Theta_{A'}=\Theta_A$, and a commutative square
        \begin{equation*}
        \xymatrix{
        E \ar[r]^h \ar[d]^{f_E} & A \ar[d]^{f_A} \\
        E' \ar[r]^{h'} & A'.
        }
        \end{equation*}

The moduli space $\wNL^b_{g,d}(N)$ may be constructed from $\wNL_{g,d}=\wNL_{g,d}(1)$, as considered in \cite{gl}, as a union of connected components of the fiber product $\wNL_{g,d}\times_{(\cA_g\times\cA_1)} (\cA_g(N)\times\cA_1(N))$. In particular, $\wNL^b_{g,d}(N)$ is smooth of dimension $\binom{g}{2}+1$. By duality, $\wNL^b_{g,d}(N)$ may equivalently be viewed as the moduli space of maps $h^\vee:A\to E$ respecting level structure.

Let $\epsilon^b(N):\wNL^b_{g,d}(N)\to \cA_g(N)$ and $\mu^b(N):\wNL^b_{g,d}(N)\to \cA_1(N)$ be the forgetful maps remembering the target and source, respectively, of $h$.
\end{defn}

There is a map $\nu_N:\wNL^b_{g,d}(N)\to \wNL_{g,d}$ forgetting level structure, which is surjective onto a union of components of $\wNL_{g,d}$. For example, if $b$ is injective, then $\nu_N$ surjects onto components of $\wNL_{g,d}$ parametrizing $h:E\to A$ that are injective on $N$-torsion. However, if $\gcd(d,N)>1$, then there are components of $\wNL_{g,d}$ parametrizing $h:E\to A$ that factor through an isogeny $E\to E'$ of degree dividing $N$, which are thus not in the image of $\nu_N$. On the other hand, if $b=0$ and $\gcd(d,N)=1$, then $\wNL^b_{g,d}(N)$ is empty.

By the same proof as in \cite[Lemma 3.4]{gl} (the level structure does not affect the arguments), the morphism $\epsilon^b(N)$ is proper, as is the morphism $(\epsilon^b(N),\mu^b(N)):\wNL^b_{g,d}(N)\to \cA_1(N)\times \cA_g(N)$. Thus, one may consider the cycle classes associated to these morphisms.

\begin{defn}
    We define the cycle classes $[\wNL^b_{g,d}(N)]\in \CH^{g-1}(\cA_g(N))$ and $[\wNL^b_{g,d}(N)]^+\in \CH^g(\cA_1(N)\times\cA_g(N))$, as well as their images in cohomology, to be the classes associated to the morphisms $\epsilon^b(N),(\mu^b(N),\epsilon^b(N))$, respectively.
\end{defn}

In \cite{gl}, the pullbacks by the pointed Torelli map $\Tor_1:\cM^{\ct}_{g,1}\to\cA_g$ of $[\wNL_{g,d}]$ and $[\wNL_{g,d}]^+$ to $\CH^{g-1}(\cM_{g,1}^{\ct})$ and $\CH^{g}(\cM^{\ct}_{g,1}\times\cM_{1,1})$, respectively, are shown to agree with the Gromov-Witten virtual classes on the moduli spaces of stable maps $\cM_{g,1}^{\ct,q}(\cE,d)$ to the universal elliptic curve, where the superscript $q$ denotes that the stable maps $f:C\to E$ are required to send the marked point of $C$ to the origin of $E$. Identical arguments show that the classes $[\wNL^b_{g,d}(N)]$ and  $[\wNL^b_{g,d}(N)]^+$ pull back to virtual classes on the spaces of stable maps with full level-$N$ structure and whose induced maps on $N$-torsion are prescribed by $b$.

\subsection{The case $d=0$}\label{d=0_sec}

Let $\lambda_i\in \CH^i(\cA_g(N))$ denote the $i$-th Chern class of the Hodge bundle on $\cA_g(N)$. By convention, we set
\begin{align*}
    [\wNL^b_{g,0}]&=(-1)^{g}\frac{1}{24}\lambda_{g-1}\in \CH^{g-1}(\cA_g(N)),\\
    [\wNL^b_{g,0}]^+&=(-1)^g\lambda_g=0\in \CH^g(\cA_1(N)\times\cA_g(N))
\end{align*}
 if $b\equiv0\pmod{N}$, and both cycle classes $[\wNL^b_{g,0}],[\wNL^b_{g,0}]^+$ to be zero otherwise. (See \cite[Proposition 1.2]{vandergeer} for the vanishing of $\lambda_g$.)
 
 The definitions are explained as follows. The moduli space of maps $f:E\to A$ of degree zero inducing $b$ on $H_1$ is empty if $b\not\equiv0\pmod{N}$, and isomorphic to $\cA_1(N)\times\cA_g(N)$ if $b\equiv0\pmod{N}$. We also add cusps in the first factor, passing to $\cA_1(N)^{*}\times\cA_g(N)$. A constant map $f:E\to A$ has obstruction space 
 \begin{equation*}
 H^1(E,f^{*}T_A)\cong H^1(E,\cO_E)\otimes T_0A\cong H^0(E,\omega_E)^\vee \otimes H^0(A,\Omega_A)^\vee.
 \end{equation*}

 Thus, the product $\cA_1(N)^{*}\times\cA_g(N)$ is equipped naturally with the global obstruction bundle $\bE_1^{\vee}\otimes\bE_g^{\vee}$ given by the tensor product of Hodge bundles on each factor. The virtual class $[\wNL^b_{g,0}]^+\in \CH^g(\cA_1(N)^{*}\times\cA_g(N))$ is therefore naturally given by the top Chern class
 \begin{align*}
     c_g(\bE_1^{\vee}\otimes\bE_g^{\vee})&=c_g(\bE_g^{\vee})+c_1(\bE_1^{\vee})c_{g-1}(\bE_g^{\vee})\\
     &=(-1)^g\lambda_g+\left(-\frac{1}{24}\cdot[E_0]\right)\cdot((-1)^{g-1}\lambda_{g-1}),
 \end{align*}
 where $[E_0]\in \CH^1(\cA_1^{*})$ is the class of a geometric point. Pushing forward to $\CH^{g-1}(\cA_g(N))$ gives the formula for $[\wNL^b_{g,0}]\in \CH^{g-1}(\cA_g(N))$ above, and restricting to the interior gives the formula for $ [\wNL^b_{g,0}]^+\in \CH^g(\cA_1(N)\times\cA_g(N))$.

 \subsection{Vanishing}

\begin{prop}\label{N=1_trivial}
The classes $[\wNL^b_{g,d}(N)]^+\in H^{2g}(\cA_1(N)\times\cA_g(N))$ in cohomology are supported in the odd K\"{u}nneth component $H^{1}(\cA_1(N))\otimes H^{2g-1}(\cA_g(N))$. In particular, they vanish when $N=1$.
\end{prop}

\begin{proof}
    We may assume that $d\ge 1$. Passing from $\cA_1(N)$ to $\cA_1(N)^*$, we have that the map $\wNL^b_{g,d}(N)\to \cA_1(N)^*\times\cA_g(N)$ is proper, by the same argument as in \cite[Lemma 3.4]{gl}. Thus, we may consider the cycle class $[\wNL^b_{g,d}(N)]^+$ as an element of $H^{2g}(\cA_1(N)\times\cA_g(N))$ or of $H^{2g}(\cA_1(N)^*\times\cA_g(N))$; the pullback of the former is equal to the latter. 
    
    Consider now the projection of $[\wNL^b_{g,d}(N)]\in H^{2g}(\cA_1(N)^*\times\cA_g(N))$ to the K\"{u}nneth component $H^{0}(\cA_1(N)^*)\otimes H^{2g}(\cA_g(N))$. Up to a constant, the projection map is given by $\alpha\mapsto 1\otimes p_{*}(\alpha\cap ([E]\times\cA_g(N)))$, where $[E]\in H^2(\cA_1(N)^*)$ is the class of any point and $p_{*}:H^{2g+2}(\cA_1(N)^*\times\cA_g(N))\to H^{2g}(\cA_g(N))$ is proper pushforward. Indeed, this map is easily seen to be the identity on $H^{0}(\cA_1(N)^*)\otimes H^{2g}(\cA_g(N))$ and zero on the other two K\"{u}nneth components.

    In particular, we may take $[E]$ to be any cusp of $\cA_1(N)^{*}$, and in this case we have that $\wNL^b_{g,d}(N)\cap ([E]\times \cA_g(N))$ is empty in $\cA_1(N)^*\times\cA_g(N)$, because smooth PPAVs contain no rational curves. Thus, $[\wNL^b_{g,d}(N)]^+$ projects to zero in $H^{0}(\cA_1(N)^*)\otimes H^{2g}(\cA_g(N))$, so the same is true in $H^{0}(\cA_1(N))\otimes H^{2g}(\cA_g(N))$. Moreover, the K\"{u}nneth component $H^{2}(\cA_1(N))\otimes H^{2(g-1)}(\cA_g(N))$ is identically zero, so the claim follows.
    
    Finally, when $N=1$, we have $H^{1}(\cA_1)=0$, as $\cA_1$ is contractible. Thus, the class $[\wNL_{g,d}]^+
    \in H^{2g}(\cA_1\times\cA_g)$ is identically zero.
\end{proof}

\begin{rem}\label{chow_vanish_remark}
When $N=1$, the same argument shows that the classes $[\wNL_{g,d}]^+$ vanish in $\CH^{g}(\cA_1\times\cA_g)$. Indeed, because the coarse space of $\cA_1^{*}$ is isomorphic to $\bP^1$, the Chow group $\CH^{g}(\cA_1^{*}\times\cA_g)$ admits a K\"{u}nneth decomposition
\begin{equation*}
    \CH^{g}(\cA^{*}_1\times\cA_g)\cong \CH^{1}(\cA^{*}_1)\otimes \CH^g(\cA_g)\oplus  \CH^{0}(\cA^{*}_1)\otimes \CH^{g-1}(\cA_g),
\end{equation*}
and one can proceed as in the proof of Proposition \ref{N=1_trivial}.
\end{rem}

\begin{rem}\label{N<=5}
In fact, we have $[\wNL^b_{g,d}(N)]^+=0$ in both cohomology and Chow for any $N\le 5$, because $\cA_1(N)$ is rational in this range. In the appendix, it is shown that $[\wNL^b_{2,d}(N)]^+$ is not always zero when $N\ge 11$ or $N\neq12$. For the remaining values of $N$, we do not know whether the classes $[\wNL^b_{g,d}(N)]^+$ always vanish.
\end{rem}

\section{Uniformization}

Let $\Gamma(N)\subset\SL(2,\bZ)$ be the subgroup of matrices congruent to the identity modulo $N$, and more generally let $\Gamma(N)_g\subset\Sp(2g,\bZ)$ be the subgroup of such matrices for any $g$. Then, we have $\cA_1(N)=\Gamma(N)\backslash\bH$ and $\cA_g(N)=\Gamma(N)_g\backslash\bH_g$.

We also denote by
\begin{equation*}
    J_2=
    \begin{pmatrix}
        0 & 1 \\
        -1 & 0
    \end{pmatrix}
    ,
    J_{2g}=
        \begin{pmatrix}
        0 & \Id \\
        -\Id & 0
    \end{pmatrix}
\end{equation*}
the matrices of the standard symplectic forms, where $\Id$ is the $g\times g$ identity matrix.

\begin{lem}\label{borel}
Let $M_{2g\times 2,d}$ be the set of integer $2g\times 2$ matrices whose columns span a rank 2 sublattice of discriminant $d^2$ in the standard symplectic lattice $\bZ^{2g}$.

Then, for any $N\ge1$, the set $\Gamma(N)\backslash M_{2g\times 2,d}/\Gamma(N)_g$ is finite.
\end{lem}

\begin{proof}
Let $V_\bQ= \bQ\otimes_\bZ \bZ^{2g}$ be the standard symplectic $\bQ$-vector space. The symplectic group $G_\bQ=\Sp(2g,\bQ)$ acts on $V$, so it also acts diagonally on $W=V\oplus V$. Let $O$ be the $G_\bQ$-orbit of an element of $M_{2g\times 2,d}$; this orbit clearly contains all of $M_{2g\times 2,d}$. By \cite[Theorem 9.11]{borel}, $O\cap W_{\bZ}$ is composed of finitely many orbits for $G_\bZ$.
\end{proof}

\begin{lem}\label{lem:B_discriminant}
Let $E=\bC/\Lambda_2$ be an elliptic curve and let $A=\bC^g/\Lambda'_{2g}$ be a PPAV. Choose symplectic bases $H_1(E)\simeq \Lambda_{2}$ and $H_1(A)\simeq \Lambda'_{2g}$ with respect to the polarization forms.

Let $h:E\to A$ be a map of degree $d$, and let $B_h$ be the matrix of the induced map on homology $H_1(E)\to H_1(A)$ with respect to the chosen bases of $\Lambda_2,\Lambda'_{2g}$. Then, we have $B_h\in M_{2g\times 2,d}$.
\end{lem}

\begin{proof}
If $h:E\to A$ has degree $d$, then the composition 
\begin{equation*}
    E \overset{h}\longrightarrow A \to A^\vee \overset{h^\vee}\longrightarrow E^\vee \to E
\end{equation*}
is equal to the multiplication by $d$ map $[d]$. The induced maps on first homology groups are given by
\begin{equation*}
H_1(E)\to H_1(A) \overset{J_{2g}}\longrightarrow H_1(A^\vee) \to H_1(E^\vee) \overset{J_2^{-1}}\longrightarrow H_1(E),
\end{equation*}
and the composition is $d\Id$,

With respect to their chosen bases, the map $H_1(E) \to H_1(A)$ is given by $B_h\in M_{2g\times 2}$, and the map $H_1(A^\vee) \to H_1(E^\vee)$ by its transpose $B_h^T$, so we have
\begin{equation*}
J_2^{-1} B_h^T J_{2g} B_h = d\Id.
\end{equation*}
This is equivalent to
\begin{equation*}
B_h^T J_{2g} B_h = dJ_2,
\end{equation*}
so $B_h\in M_{2g\times 2, d}$.
\end{proof}

\begin{lem}\label{lem:B_vanishing}
Given a pair $(\tau,\tau')\in \bH\times \bH_g$, there exists a degree $d$ map $E_\tau \to A_{\tau'}$ between the associated abelian varieties if and only if 
\begin{equation*}
    (B \tilde{\tau})^T\cdot \tilde{\tau}'=0\in\bC^g
\end{equation*}
for some matrix $B\in M_{2g\times 2,d}$. Here, $\tilde{\tau} = \begin{pmatrix} \tau \\ \Id \end{pmatrix}$ denotes the Siegel augmentation.
\end{lem}
\begin{proof}
Let $h: E \to A$ be a degree $d$ map. Let $B_h\in M_{2g\times 2,d}$ be the matrix for 
$$h_*:H_1(E) \to H_1(A)$$
with respect to symplectic bases $\langle\alpha,\beta\rangle$ for $H_1(E)$ and $\langle\alpha_j,\beta_j\rangle_{j=1,\ldots,g}$ for $H_1(A)$. The graph $\Gamma_h \subset E \times A$ has homology class
\begin{equation*}
\alpha\times h_*\beta - \beta\times h_*\alpha + E\times [\pt] + [\pt]\times h_*[E]\in H_2(E\times A,\bZ).
\end{equation*}

If $\omega\in H^{1,0}(E)$ and $\omega_j\in H^{1,0}(A)$ are the normalized holomorphic 1-forms, then their external wedge product $\omega \wedge \omega_j\in H^{2,0}(E\times A)$ integrates to 0 on any algebraic curve, so on $\Gamma_h$ in particular. In terms of the decomposition above, this implies that
\begin{equation*}
\int_\alpha \omega \int_\beta h^*\omega_j - \int_\beta \omega \int_\alpha h^*\omega_j=0
\end{equation*}
for $j=1,2,\dots,g$. 

The Siegel augmented period matrices are given by
\begin{align*}
    \tilde{\tau} &= \begin{pmatrix} \int_\alpha \omega \\ \int_\beta \omega  \end{pmatrix},\\
    \tilde{\tau}' &= \begin{pmatrix} \int_{\alpha_i} \omega_j \\ \int_{\beta_i} \omega_j  \end{pmatrix}.
\end{align*}
Multiplying $\tilde{\tau}'$ by $B_h^T$ on the left has the effect of replacing $\alpha_i$ and $\beta_j$ with the pushforwards of $\alpha$ and $\beta$ under $h$:
\begin{equation*}
B_h^T \tilde{\tau}' = 
\begin{pmatrix} 
\int_{h_*\alpha} \omega_j \\ 
\int_{h_*\beta} \omega_j  \end{pmatrix} 
=
\begin{pmatrix} 
\int_{\alpha} h^*\omega_j \\ 
\int_{\beta} h^*\omega_j 
\end{pmatrix}
.
\end{equation*}
Now, we also have
\begin{equation*}
J_2^{-1}\tilde{\tau} 
=
\begin{pmatrix}
- \int_\beta \omega \\ 
\int_\alpha \omega 
\end{pmatrix}
\end{equation*}
and the vanishing above is equivalent to 
\begin{equation*}
    \tilde{\tau}^T J_2 B_h^T \tilde{\tau}'= (B_hJ_2^{-1}\tilde{\tau})^T\cdot \tilde{\tau'} =0.
\end{equation*}
Taking $B=B_hJ_2^{-1}$ yields the first direction of the Lemma.

Conversely, given $B\in M_{2g\times 2,d}$, define $B_h=BJ_2$, and a linear subtorus 
\begin{equation*}
    \Gamma\subset E\times A=\bC/\Lambda_\tau\times \bC^g/\Lambda'_{\tau'}
\end{equation*}
by $\Gamma=\{(z,B_hz)\mid z\in\bC\}$, where we extend the symplectic bases of the lattices $\Lambda_\tau,\Lambda'_{\tau'}$ to $\bC,\bC^g$, respectively. Reversing the previous calculation, the vanishing $(B \tilde{\tau})^T\cdot \tilde{\tau}'=0$ implies that integrals of holomorphic 2-forms on $\Gamma$ vanish, which in turn implies that $\Gamma\cong E$ is a complex subtorus. Post-composing with the projection to $A$ gives the desired map $h:E\to A$.
\end{proof}

\begin{cor}
    We have an isomorphism
    \begin{equation*}
        \wNL^b_{g,d}(N)\cong \Gamma(N)\backslash \left\{(\tau,\tau',B)\in\bH\times\bH_g\times M_{2g\times 2,d} \middle\vert \, 
        \begin{array}{@{}c@{}}
            (B \tilde{\tau})^T\cdot \tilde{\tau}'=0\\
            b\equiv BJ_2\pmod{N}
        \end{array}
        \right\} /\Gamma(N)_g
    \end{equation*} 
\end{cor}

\begin{proof}
Lemmas \ref{lem:B_discriminant} and \ref{lem:B_vanishing} show that, given $E_\tau\in \cM_{1,1}(N)$ and $A_{\tau'}\in \cA_g(N)$, the data of $f:E_\tau\to A_{\tau'}$ of degree $d$ is equivalent to the data of a matrix $B\in M_{2g\times 2,d}$, up to the actions of $\Gamma(N),\Gamma(N)_g$, satisfying $(B \tilde{\tau})^T\cdot \tilde{\tau}'=0$. Moreover, the induced map $b:(\bZ/N\bZ)^{2}\to(\bZ/N\bZ)^{2g}$ on $N$-torsion is given, by the calculation of Lemma \ref{lem:B_vanishing}, by the matrix $B_h=BJ_2$. 
\end{proof}

\section{Kudla-Millson Modularity}\label{km_sec}
\subsection{Symmetric spaces}
Let $(W_2, \omega)$ and $(W'_{2g}, \omega')$ be real symplectic vector spaces. The tensor product $V=W_2\otimes W'_{2g}$ has a natural symmetric pairing $\gamma$ given by the product $\omega\cdot\omega'$. Note that all pure tensors in $V$ are isotropic with respect to $\gamma$. Choose Darboux bases $
\langle e,f\rangle$ and $\langle e_i,f_i\rangle$ for $W_2$ and $W'_{2g}$, respectively, so that
$$e\otimes e_i' + f\otimes f_i'\,\,\,(1\leq i \leq g)$$
$$e\otimes f_i' - f\otimes e_i'\,\,\,(1\leq i \leq g)$$
form a basis for a maximal positive definite subspace $P_0\subset V$. Similarly,
$$e\otimes e_i' - f\otimes f_i'\,\,\,(1\leq i \leq g)$$
$$e\otimes f_i' + f\otimes e_i'\,\,\,(1\leq i \leq g)$$
form a basis for a maximal negative definite subspace $N_0\subset V$, with $N_0 = P_0^\perp$. Hence, the symmetric pairing $\gamma$ is non-degenerate of signature $(2g,2g)$.

Next, consider the map 
\begin{equation*}
\varphi:\SL_2(\bR) \times \Sp_{2g}(\bR) \to SO(V)_0 \simeq SO(2g,2g)_0
\end{equation*}
defined by $\varphi(M,M')=M\otimes M'$. Its kernel is $\{\pm (\Id,\Id)\}$, which is contained in the maximal compact $K=SO(2)\times U(g)$. The restriction of $\varphi$ to $K$ lands in $K'$:
$$\varphi|_K: K \to K'=SO(2g)\times SO(2g)\subset SO(2g,2g)_0.$$
Hence $\varphi$ induces an embedding $\phi$ on the associated symmetric spaces:
$$\phi: \bH \times \bH_g \to \Gr^-(2g,V).$$

The symmetric space for $SO(2g,2g)_0$ may be identified with the positive definite Grassmannian or the negative definite Grassmannian; we choose the latter. Explicitly, given $(\tau,\tau')\in \bH \times \bH_g$, there exist matrices $(M,M')\in \SL_2(\bR) \times \Sp_{2g}(\bR)$ sending $(i,iI_g)\mapsto (\tau,\tau')$.
$$\phi(\tau,\tau') = (M\otimes M')(N_0)\in \Gr^-(2g,V).$$
One easily checks that $\Stab(i,iI_g)$ preserves $N_0$, so the map $\phi$ is well-defined.\\

\begin{prop}\label{siegelcomp}
Let 
\begin{equation*}
B=\begin{pmatrix}
a_1 & b_1 \\
a_2 & b_2 \\
\vdots & \vdots\\
a_{2g} & b_{2g}
\end{pmatrix}
\in M_{2g\times 2}
\end{equation*}
be an integer $2g\times 2$ matrix. Let
\begin{equation*}
    B_{\phi}=\sum_{k=1}^g b_{g+k}(e\otimes e'_k)-b_k(e\otimes f'_k)-a_{g+k}(f\otimes e'_k)+a_k(f\otimes f'_k)\in W_2\otimes W'_{2g}.
\end{equation*}

Then, for any $(\tau,\tau')\in \bH\times \bH_g$, the following are equivalent:
\begin{itemize}
\item $\phi(\tau,\tau')$ is orthogonal to $B_\phi$ in $V$.
\item $(B\tilde{\tau})^T\cdot \tilde{\tau}'=0\in \bC^g$. Here,  $\tilde{\tau} = \begin{pmatrix} \tau \\ \Id \end{pmatrix}$ denotes the Siegel augmentation.
\end{itemize}
\end{prop}

\begin{proof}
Let
\begin{equation*}
    M=\begin{pmatrix}w & x\\ y & z\end{pmatrix},
    M'=\begin{pmatrix}W & X\\ Y& Z\end{pmatrix},
\end{equation*}
 where $W,X,Y,Z$ are $g\times g$ real matrices. then, we have
 \begin{equation*}
     \tau=\frac{wi+x}{yi+z},\tau'=(iW+X)(iY+Z)^{-1}.
 \end{equation*}


For $j=1,2,\ldots,g$, the $j$-th entry of the $1\times g$ matrix $(yi+z)(B\tilde{\tau})^T \cdot \tilde{\tau}'(iY+Z)$
is
\begin{equation*}
    \beta_j=\sum_{k=1}^{g}(a_k(wi+x)+b_k(yi+z))(iW+X)_{kj}+\sum_{k=1}^{g}(a_{g+k}(wi+x)+b_{g+k}(yi+z))(iY+Z)_{kj}
\end{equation*}
where $(iW+X)_{kj},(iY+Z)_{kj}$ denote the entries in the $k$-th row and $j$-th column of the respective matrices. Thus, we have
\begin{align*}
    \Re(\beta_j)&=\sum_{k=1}^g [a_k(xx_{kj}-ww_{kj})+b_k(zx_{kj}-yw_{kj})+a_{g+k}(xz_{kj}-wy_{kj})+b_{g+k}(zz_{kj}-yy_{kj})],\\
    \Im(\beta_j)&=\sum_{k=1}^g [a_k(wx_{kj}+xw_{kj})+b_k(yx_{kj}+zw_{kj})+a_{g+k}(wz_{kj}+xy_{kj})+b_{g+k}(yz_{kj}+zy_{kj})].
\end{align*}

On the other hand, we compute that $\phi(\tau,\tau')=(M\otimes M')(N_0)$ is spanned by
\begin{align*}
    r_j&=(we+yf)\otimes(We'_j+Yf'_j)-(xe+zf)\otimes(Xe'_j+Zf'_j),\\
    &=-\sum_{k=1}^g [(e\otimes e'_k)(xx_{kj}-ww_{kj})+(e\otimes f'_k)(xz_{kj}-wy_{kj})\\
    &\quad\quad+(f\otimes e'_k)(zx_{kj}-yw_{kj})+(f\otimes f'_k)(zz_{kj}-yy_{kj})],\\
    s_j&=(we+yf)\otimes(Xe'_j+Zf'_j)+(xe+zf)\otimes(We'_j+Yf'_j),\\
    &=\sum_{k=1}^g [(e\otimes e'_k)(wx_{kj}+xw_{kj})+(e\otimes f'_k)(wz_{kj}+xy_{kj})\\
    &\quad\quad+(f\otimes e'_k)(yx_{kj}+zw_{kj})+(f\otimes f'_k)(yz_{kj}+zy_{kj})].
\end{align*}
for $j=1,2,\ldots,g$. Here, the matrices $W,X$ are taken to act on the basis $\{e'_1,\ldots,e'_g\}$ and the matrices $Y,Z$ are taken to act on the basis $\{f'_1,\ldots,f'_g\}$.

Finally, we see that $\gamma(B_\phi,-r_j)=\Re(\beta_j)$ and $\gamma(B_\phi,s_j)=\Im(\beta_j)$, which yields the needed equivalence.
\end{proof}

\subsection{Arithmetic Quotients}
Fix unimodular lattices $\Lambda\subset W_2$ and $\Lambda' \subset W'_{2g}$, which give rise to an even unimodular lattice $L:=\Lambda\otimes \Lambda' \subset V$. The map of $\phi$ between symmetric spaces descends to a map on arithmetic quotients: 
\begin{equation*}
\phi:\Aut(\Lambda)_N \backslash \bH \times \Aut(\Lambda')_N\backslash \bH_g \to \Aut(L)_N\backslash \Gr^-(2g,V)=:\km(L(N)).
\end{equation*}
Here, $\Aut(\Lambda)_N,\Aut(\Lambda')_N,\Aut(L)_N$ denote the automorphism groups of the respective lattices that reduce to the identity mod $N$, so that $\phi$ is a map $\cA_1(N)\times\cA_g(N)\to\km(L(N))$.

If $v\in L$ is a vector of positive norm, then set
\begin{equation*}
\wt{C}_v := \{P\in \Gr^-(2g,V): P\subset v^\perp\}.
\end{equation*}
This is a non-empty sub-symmetric space of codimension $2g$. By Lemma \ref{borel}, $\Aut(L)_N$ acts on the lattice vectors of norm $2d>0$ with finitely many orbits. 

Fix now an abelian group homomorphism $b:(\bZ/N\bZ)^2\to (\bZ/N\bZ)^{2g}$. Choose symplectic bases $\{e,f\}$ and $\{e'_1,\ldots,e'_g,f'_1,\ldots,f'_g\}$ of $\Lambda,\Lambda'$, respectively, and given
\begin{equation*}
    v=\sum_{k=1}^g b_{g+k}(e\otimes e'_k)-b_k(e\otimes f'_k)-a_{g+k}(f\otimes e'_k)+a_k(f\otimes f'_k)\in L
\end{equation*}
define
\begin{equation*}
    v_h=
    \begin{pmatrix}
        -b_1 & a_1\\
        -b_2 & a_2\\
        \vdots & \vdots\\
        -b_{2g} & a_{2g}
    \end{pmatrix}
    =\begin{pmatrix}
a_1 & b_1 \\
a_2 & b_2 \\
\vdots & \vdots\\
a_{2g} & b_{2g}
\end{pmatrix}
\cdot J_2.
\end{equation*}
Here, the vector $v\in L$ plays the role of $B_\phi$ in Propostion \ref{siegelcomp}, and $v_h$ plays the role of the matrix $B_h$ in Lemmas \ref{lem:B_discriminant} and \ref{lem:B_vanishing}. The intermediate matrix $B$ appearing in Lemma \ref{lem:B_discriminant} and Propostion \ref{siegelcomp} satisfies $BJ_2=B_h$. Given a $d$-elliptic map $h:E\to A$, the matrix of the induced map on first homology is given by $B_h$, which will be required to reduce modulo $N$ to the prescribed map $b$.

\begin{defn}
We define the special cycles
\begin{equation*}
C^b_{d}(N) := \Aut(L)_N \backslash\left(\bigcup_{\substack{v\in L /\Aut(L)_N. \\ v^2=2d \\ v_h\equiv b\pmod{N}}} \wt{C}_v \right) \subset \km(L(N)).
\end{equation*}
\end{defn}

\begin{prop}\label{special_cycle_pullback}
We have a commutative diagram
\begin{equation*}
\xymatrix{
    \wNL_{g,d}^b(N) \ar[r] \ar[d]_{(\epsilon(N),\mu(N)} & C^b_d(N) \ar@{^{(}->}[d] \\
    \cA_1(N) \times \cA_g(N) \ar[r]^(0.61){\phi} & \km(L(N))
}
\end{equation*}
inducing a birational map $\wNL_{g,d}^b(N) \to \phi^{-1}(C^b_d(N))$. In particular, we have 
\begin{equation*}
    [\wNL_{g,d}^{b}(N)]^+=\phi^{*}[ C^b_d(N)]
\end{equation*}
in $H^{2g}(\cA_1(N)\times\cA_g(N))$.
\end{prop}

\begin{proof}
Propostion \ref{siegelcomp} shows that the composition $\phi\circ(\epsilon(N),\mu(N))$ has image equal to $C^b_d(N)$, and furthermore that $\phi^{-1}(C^b_d(N))$ is equal to the locus of $(E,A)$ (with full level-$N$ structure) for which there \emph{exists} a map $h:E\to A$ inducing the map $b$ on $N$-torsion. To identify this pullback generically with $\wNL_{g,d}^b(N)$, it therefore suffices to show that the map $(\epsilon(N),\mu(N))$ is birational onto its image. 

Note that $(\epsilon(N),\mu(N))$ is unramified by \cite[Proposition 3.6]{gl} (the level structure does not affect the local arguments), so it suffices to show to show that $(\epsilon(N),\mu(N))$ is generically of degree 1 on geometric points. This amounts to the statement that, at a general point $h:E\to A$ of $\wNL_{g,d}^b(N)$, the map $h$ is the only one (up to isomorphism) from $E$ to $A$. By a dimension count, we may assume that $A$ splits up to isogeny into the product of $E$ and a simple abelian variety $B$ of dimension $g-1$, and when $g=2$, we can assume further that $B$ is not isogeneous to $E$. We may furthermore assume that $E$ is general. 

Now, if there are two non-isomorphic maps $h_1,h_2:E\to A$, then the induced map $E\times E\to A$ must have a 1-dimensional kernel $E'$, by the assumption on the splitting of $A$. The two maps $i_1,i_2:E'\to E$ must have the same degree $\delta$, as the two compositions $E'\to E\times\{0\} \to A$ and $E'\to \{0\}\times E\to A$ agree up to sign, and $\deg(h_1)=\deg(h_2)=d$. Thus, the two maps $i_1\circ i_1^\vee$ and $i_2\circ i_1^\vee$ must both have map $\delta^2$. The first map is multiplication by $\delta$. If $E$ is general, then $\text{End}(E)\cong\bZ$, so $i_2\circ i_1^{\vee}$ must either be multiplication by $\delta$ or $-\delta$. In particular, we have $i_1\circ i_1^{\vee}=\pm i_2\circ i_2^\vee$, and hence $i_1=\pm i_2$, implying that $h_1,h_2$ are the same geometric point of $\wNL_{g,d}^b(N)$.
\end{proof}

The main input into Theorem \ref{thm_levelN} is the modularity of Kudla-Millson.

\begin{thm}[Kudla-Millson \cite{km}]
Let $e_0$ be the Euler class of the dual of the rank $2g$ tautological bundle on $\km(L(N))$ of negative definite $2g$-planes. For any $b\in L(N)^\vee/L(N)$, the power series
\begin{equation*}
    \Phi^{b}(q):=e_0+\sum_{d\ge1}[C^b_d(N)]q^d \in \Mod(2g,\Gamma(N))\otimes H^{2g}(\km(L(N)))
\end{equation*}
is a modular form of weight $2g$ and level $N$.
\end{thm}

\begin{proof}[Proof of Theorem \ref{thm_levelN}]
The pullback under $\phi$ of the tautological bundle is a real vector bundle on $\A_1\times \A_g$ whose complexification has fiber at $(E,X)$ equal to $H^{1,0}(E)\otimes H^{1,0}(X) \oplus H^{0,1}(E)\otimes H^{0,1}(X)$.  The two direct summands are duals and complex conjugates of each other. The inclusion of the real $2g$-plane followed by projection onto the first summand gives an isomorphism of real oriented vector bundles, so $$\phi^*(e_0) = (-1)^g c_g(\bE_1 \boxtimes \bE_g) \in H^{2g}(\A_1\times \A_g,\bQ).$$
\end{proof}

\appendix \label{appendix}
\section{Nonvanishing of certain Noether-Lefschetz classes}
\begin{center}
by N. Sweeting
\end{center}

\subsection{Overview}
In this appendix, we prove the nonvanishing of the Noether-Loefschetz classes $[\widetilde{\mathsf{NL}}_{2,1}^b(N)]^+$ for sufficiently large $N$ when $b$ is an embedding (Theorem \ref{main appendix theorem} below).
The strategy of the proof is to produce, using theta lifts from $\operatorname{GSO}_{2,2}$ to $\GSp_4$, explicit classes in $H^4_c(\mathcal A_1(N)\times \mathcal A_2(N), \C)$ with nonzero pairing against $[\widetilde{\mathsf{NL}}_{2,1}^b(N)]^+$ under Poincar\'e duality. Most of the theoretical work is contained in \cite[Theorems A and C]{sweeting2022tateclassesendoscopyoperatornamegsp4}, but without any precise control of the level $N$; thus we must supplement the methods of \emph{loc. cit.} with a number of additional computations to make the level structures explicit. 

\subsection{Conventions}
\subsubsection{} If $G$ is an algebraic group over $\Q$, then $[G] \colonequals G(\Q)\backslash G(\bA)$ denotes the usual adelic quotient. We denote by $\mathcal A(G)$ the space of automorphic forms on $[G]$, and by $\mathcal A_0(G)$  the subspace of cusp forms. If $K \subset G(\bA_f)$ is a compact open subgroup, then we write $\mathcal A_0(G; K)$ for the space of $K$-invariant cusp forms.
\subsubsection{}
For an integer $N\geq 1$, we consider the compact open subgroup $$K_1(N) = \prod K_1(N)_p= \set{\begin{pmatrix}
    a & b \\ c & d 
\end{pmatrix}\in \GL_2(\widehat \Z)\,:\, c \in N\widehat\Z ,\; d\in 1 + N \widehat\Z }$$
of $\GL_2(\bA_f)$. It is clear that $K_1(N)_p$ depends only on the $p$-adic valuation of $N$.
\subsubsection{}
We denote by $B \subset \GL_2$ the upper triangular Borel subgroup, and by $U \subset B$ the unipotent radical. We define a map of algebraic groups $\mathbb G_m \to \GL_2$ by $c\mapsto h_c =\begin{pmatrix}  1 & 0 \\ 0 & c \end{pmatrix}.$
\subsubsection{}
Let $\psi: \mathbb Q \backslash\mathbb A \to \C$ be the unique everywhere unramified character such that $\psi(x) = e^{2\pi i x}$ for $x\in \R$, and let $\psi_k $ be the local component of $\psi$ for every completion $k$ of $\Q$. 
\subsubsection{}
Let $\SO(2) \subset\GL_2(\R)$ be the standard maximal compact subgroup; we denote by $\chi_m$ the weight-$m$ character of $\SO(2)$ defined by $$\begin{pmatrix}
    \cos \theta & \sin\theta \\ - \sin\theta & \cos\theta
\end{pmatrix}\mapsto (\cos\theta + i \sin\theta)^m.$$
\subsection{Shimura varieties}\label{appendix subsection shimura varieties}
\subsubsection{Disconnected moduli spaces of abelian varieties}
Fix $N\geq 1$, and define $K_N = K_{N,g}= \prod_p K_{N,g,p} \subset \GSp_{2g}(\widehat \Z)$ to be the compact open subgroup of matrices that are congruent to the identity modulo $N$.  Let $\mathcal A_g'(N)$ be the complex Shimura variety for $\GSp_{2g}$ of level $K_N$:
\begin{equation}
    \mathcal A_g'(N)= \GSp_{2g}(\Q) \backslash \GSp_{2g} (\bA_f) \times \mathbb H_g/K_N.
\end{equation}
There is a natural projection $\mathcal A_g'(N) \to \Q^\times \backslash \bA_f^\times/ (1 + N \widehat \Z)^\times\simeq \mu_N$, whose fibers are the geometric connected components, each isomorphic to $\mathcal A_g(N)$. 
We also have the natural embedding
\begin{equation}
    \mathcal A_{1}'(N) \times_{\mu_N} \mathcal A_{g-1}'(N) \hookrightarrow \mathcal A_1'(N)\times \mathcal A_g'(N)
\end{equation}
corresponding to the embedding of groups $\GSp_2 \times_{\mathbb G_m} \GSp_{2g-2} \hookrightarrow\GSp_{2}\times \GSp_{2g}.$

\begin{prop}\label{appendix prop:rewrite NL class}
    For all $b: (\Z/N\Z)^2 \hookrightarrow (\Z/N\Z)^{2g}$, we have 
    $$0 \neq [\widetilde{\mathsf{NL}}_{g,1}^b(N)]^+ \in H^4(\mathcal A_1(N)\times \mathcal A_g(N), \Q)$$ if and only if $$0 \neq [\mathcal A_{1}'(N) \times_{\mu_N} \mathcal A_{g-1}'(N)] \in H^4(\mathcal A_1'(N)\times \mathcal A_g'(N), \Q).$$
\end{prop}
\begin{proof}
    Because all choices of the embedding $b$ differ by an element of $\Sp_4(\Z/N\Z)$, the resulting classes $[\widetilde{\mathsf{NL}}_{g,1}^b]^+$ are transitively permuted by the natural action of $\Sp_4(\Z/N\Z)$ on $\mathcal A_g(N)$. Hence, we may assume without loss of generality that $b$ is any fixed embedding. 
    
    The embedding $\mathcal A'_{1}(N)\times_{\mu_N} \mathcal A'_{g-1}(N)\hookrightarrow\mathcal A'_1(N)\times \mathcal A'_g(N)$ factors through the open and closed subvariety $\mathcal A_1'(N)\times _{\mu_N}\mathcal A_g'(N)\subset \mathcal A_1'(N)\times \mathcal A_g'(N)$. For each connected component $\mathcal A_1(N) \times \mathcal A_g(N)$ of $\mathcal A_1'(N)\times _{\mu_N}\mathcal A_g'(N)$, with a good choice of $b$ we have that $$(\mathcal A_{1}'(N) \times_{\mu_N} \mathcal A_{g-1}'(N))\cap (\mathcal A_1(N) \times \mathcal A_g(N)) = \widetilde{\mathsf{NL}}_{g,1}^b,$$ and the proposition follows. 
    
\end{proof}
\begin{lem}\label{appendix lemma:divisible ok}
    Let $N, M\geq 1$ be integers, and suppose $$0 \neq  \left[\mathcal A_1'(N)\times_{\mu_N} \mathcal A_{g-1}'(N)\right] \in H^4(\mathcal A_1'(N) \times \mathcal A_g'(N),\Q).$$ Then $$0 \neq  \left[\mathcal A_1'(NM)\times_{\mu_{NM}} \mathcal A_{g-1}'(NM)\right] \in H^4(\mathcal A_1'(NM) \times \mathcal A_g'(NM),\Q).$$
\end{lem}
\begin{proof}
    The pullback map $H^4(\mathcal A_1'(N) \times \mathcal A_g'(N),\Q)\to H^4(\mathcal A_1'(NM) \times \mathcal A_g'(NM),\Q)$ is injective because the projection $\pi_{N,M}:\mathcal A_1'(NM) \times \mathcal A_g'(NM)\to \mathcal A_1'(N) \times \mathcal A_g'(N)$ is finite. Hence by assumption, $\pi_{N,M}^\ast\left[\mathcal A_1'(N)\times_{\mu_N} \mathcal A_{g-1}'(N)\right] \neq 0$. On the other hand, the preimage $\pi_{N,M}^{-1} (\mathcal A_1'(N)\times_{\mu_N} \mathcal A_{g-1}'(N))$ is a union of $\GSp_2(\Z/NM\Z) \times \GSp_{2g}(\Z/NM\Z)$-translates of $\mathcal A_1'(NM)\times_{\mu_{NM}} \mathcal A_{g-1}'(NM)$, and the lemma follows. 
\end{proof}
\subsubsection{Automorphic forms in cohomology}\label{appendix subsec:automorphic forms}
Fix $\tau=  i \operatorname{Id} \in \mathbb H_g$. The stabilizer of $\tau$ in $\GSp_{2g}(\R)$ is the subgroup $\R^\times \cdot U(g)$, and the tangent space to the real manifold $\mathbb H_g$ at $\tau$ is $\mathfrak p \colonequals \mathfrak{sp}_{2g, \R}/ \mathfrak{u}(g)$. As a $U(g)$-module, $\mathfrak p_\C$ is isomorphic to $\Sym^2 \oplus (\Sym^2)^\vee$, where $\Sym^2$ is the symmetric square of the $g$-dimensional defining representation of $U(g)$. 

Thus we have a canonical map
\begin{equation}
    (\mathcal A(\GSp_{2g})\otimes \wedge^i \mathfrak p^\ast_\C)^{\R^\times \cdot U(g)} \to H^i(\mathcal A_g'(N), \C),
\end{equation}
which on cusp forms restricts to a map 
\begin{equation}\label{appendix eqn:cusp form realization map}
    (\mathcal A_0(\GSp_{2g})\otimes \wedge^i \mathfrak p^\ast_\C)^{\R^\times \cdot U(g)} \to H^i_c(\mathcal A_g'(N), \C).
\end{equation}
\subsection{Newforms and Whittaker models for $\GL_2$}\label{appendix subsection some rep theory}
\subsubsection{}
Let $\pi$ be an irreducible, admissible, infinite-dimensional representation of $\GL_2(\Q_p)$ for prime $p$. 
The conductor of $\pi$ is the least $n$ such that $\pi^{K_1(p^n)_p} \neq 0$; it is well-known that such an $n$ always exists, and, if $n$ is minimal, then $\pi^{K_1(p^n)_p}$ is one-dimensional. A generator of this space is called a local newform for $\pi$. 
\subsubsection{}
Recall the nontrivial additive character $\psi_{\Q_p}$ of $\Q_p$, which is trivial on $\Z_p$ but not on $p \Z_p$. We also view  $\psi_{\Q_p}$ as a character of $$U(\Q_p)= \set{\begin{pmatrix}    1 & a \\ 0 & 1 \end{pmatrix}, \; a \in \Q_p} \subset \GL_2(\Q_p).$$
Then $\pi$ has a Whittaker model $$W_{\psi_{\Q_p}}(\pi) \subset \Ind_{U(\Q_p)}^{\GL_2(\Q_p)} \psi_{\Q_p}.$$
Let $W_{\pi,\psi_{\Q_p}}^0\in W_{\psi_{\Q_p}}(\pi)$ be a local newform. 
\begin{prop}\label{appendix prop:whittaker newform on torus}
   Suppose $\pi$ has conductor $n \geq 1$, so that $L(s,\pi) = (1 - \alpha p^{-s})$ for some $\alpha\in \C$.  Then  up to rescaling $W_{\pi,\psi_{\Q_p}}^0$, we have 
    $$W_{\pi,\psi_{\Q_p}}^0\left(\begin{pmatrix}t & 0 \\ 0 & 1\end{pmatrix} \right)= \begin{cases} 
    0, &  \ord_p(t) < 0,\\
    1, & \ord_p (t) = 0,\\
        |t|^{1/2} \alpha ^{\ord_p (t)}, &  \ord_p(t) > 0. 
    \end{cases}$$
\end{prop}
\begin{proof}
This is a special case of \cite[Theorem 4.1]{miyauchi2014whittaker}.
\end{proof}
\subsection{Induced representations and Eisenstein series on $\GL_2$}
\subsubsection{}
For each place $v$ of $\Q$,
let \begin{equation}I_v (s) = \Ind_{B(\Q_v)}^{\GL_2(\Q_v)}\delta_B^s
\end{equation}
be the normalized induction, and let $I(s) = \otimes_v 'I_v(s)$.

For $\varphi(s) \in I(s)$ a standard section, we have the Eisenstein series
$$E(g, s;\varphi) = \sum_{\gamma\in B(\Q)\backslash\GL_2(\Q)} \varphi(s) (\gamma g), g\in \GL_2(\bA),$$ which converges for $\Re(s) \gg 0$.
\subsubsection{}
For $N \geq 1$, we define a section $\varphi^0_N = \otimes_v \varphi^0_{N, v}\in I(1/2)$  as follows:
\begin{itemize}
    \item For $v = p $, $\varphi^0_{N,p}$ is the unique $K_1(N)_p$-invariant section supported on $B(\Q_p)\cdot K_1(N)_p$ and satisfying $\varphi^0_{N,p} (1) = 1.$
    \item For $v = \infty$, $\varphi^0_{N, \infty}$ is the unique $\SO(2)$-spherical section satisfying $\varphi^0_{N, \infty} (1) = 1$. 
\end{itemize}
We can extend $\varphi^0_N$ uniquely to a section $\varphi_N^0(s)\in I(s)$ so that the restriction of $\varphi^0_N$ to $\GL_2(\widehat \Z) \cdot \SO(2)$ is independent of $s$.
\begin{prop}\label{appendix prop:nonzero residue}
    The Eisenstein series $E(g, s; \varphi^0_N)$ has a pole at $s = 1/2$, with residue a nonzero constant function of $g$.
\end{prop}
\begin{proof}
    By the well-known theory of Eisenstein series for $\GL_2$, it suffices to show that $\varphi^0_{N,v}$ has nontrivial image under the intertwining operator $$M_v: I_v(1/2) \to I_v(-1/2)$$ for all primes $v$. At $v = \infty$ and $v = p\nmid N$, this is clear because $\varphi^0_{N,v}$ is the unique spherical vector for a maximal compact subgroup of $\GL_2(\Q_v)$, so we consider the case of $v = p |N$. The intertwining operator is given explicitly 
by $$M_p(\varphi)(g) = \int_{\Q_p} \varphi\left(\begin{pmatrix}
     0 & 1 \\ -1 & 0 
\end{pmatrix}\begin{pmatrix}
    1 & y \\ 0 & 1
\end{pmatrix}g\right) \d y,\;\; \varphi\in I_v(1/2), g\in \GL_2(\Q_p).$$
Normalizing the measure so that $\Z_p$ has unit volume, we therefore calculate:
\begin{align*}
    M_p(\varphi^0_{N,p})(1) &=  \varphi^0_{N,p}\left(\begin{pmatrix}
     0 & 1 \\ -1 & 0 
\end{pmatrix}\right) + \int_{\Q_p\setminus \Z_p} 
    \varphi^0_{N,p}\left(\begin{pmatrix}
     y^{-1} & 0 \\ 0 & y\end{pmatrix}
     \begin{pmatrix}
         1 & -y \\ 0 &1
     \end{pmatrix}
     \begin{pmatrix}
         1 & 0 \\ y^{-1} & 1
     \end{pmatrix}\right) \d y \\
&= 0 + \int_{\Q_p\setminus \Z_p} 
   |y|^{-2} \varphi^0_{N,p}\left(\begin{pmatrix}
         1 & 0 \\ y^{-1} & 1
     \end{pmatrix}\right) \d y \\
     &= \sum_{n\geq \ord_p(N)} \int_{p^{-n} \Z_p^\times} |y|^{-2} \d y = 
     \sum_{n \geq \ord_p(N)} p^{-n-1} (p-1)\neq 0,
\end{align*}
where in the second line we have used that $\begin{pmatrix}
    0 & 1 \\ -1 & 0
\end{pmatrix}\not\in B(\Q_p) K_1(N)_p$
for all $p|N$. 
\end{proof}
\subsection{Weil representation and theta lifting}\label{appendix subsection theta lifting}
\subsubsection{}\label{appendix bilinear forms subsubsec}
Let $\epsilon =\pm 1, $ and let $V$, $W$ be vector spaces over  a field $k$ equipped with nondegenerate $\epsilon$-symmetric and $(-\epsilon)$-symmetric pairings, respectively. We assume $\dim W = 2n $ and $\dim V = 2m $ are even, that $V$ has trivial discriminant character, and that $W$ is equipped with a complete polarization \begin{equation}
    W = W_1\oplus W_2,\;\; W_2 = W_1 ^\ast.
\end{equation}
\subsubsection{}
Let  $G_1 = G_1 (V) $, $G=G (V) $ be the connected isometry and similitude groups, respectively, of $V$, and likewise $H_1 = H _1(W) $ and $H = H (W) $; we have the natural similitude characters $\nu_G: G \to \mathbb G_m$ and $\nu_H: H \to \mathbb G_m$. 
\subsubsection{The local Weil representation}
Suppose that $k$ is a local field, and let $\psi_k$ be a nontrivial additive character of $k$. Then following Roberts' construction \cite{roberts1996theta}, the similitude Weil representation $\omega = \omega_{V, W, \psi_k}$ of $(H \times_{\mathbb G_m} G)(k)$ is realized on the Schwartz space $\mathcal  S(W_2\otimes V)$ of compactly supported, complex-valued functions on $W_2\otimes V$. Concise descriptions of this representation can be found in   \cite[\S2]{gan2011theta} or \cite[\S4]{sweeting2022tateclassesendoscopyoperatornamegsp4}, but all we require are the following two facts:
\begin{itemize}
    \item The action of $$\left(\begin{pmatrix} 1 & 0 \\ 0 & \nu_G(g) \end{pmatrix}, g\right) \in (H \times_{\mathbb G_m} G)(k) $$ on $\mathcal  S(W_2\otimes V)$ is given by $\phi \mapsto |\nu_G(g)|^{-mn/2}\phi \circ g^{-1}.$
    \item Suppose $V = V_1 \oplus V_2$ is a polarization of $V$. Then the Fourier transform
     \begin{equation}\label{appendix eqn:Fourier_transform}\begin{split}
\mathcal  S (W_2\otimes V) &\to \mathcal S(W \otimes V_2)\\
\phi&\mapsto\widehat\phi,\;\;\widehat\phi (x_1, x_2) = \int_{W_2\otimes V_1}\phi (z, x_2)\psi ( z\cdot x_1) \d z,
\end{split}
 \end{equation}
 with $\d z$ the self-dual Haar measure,
 defines an $(H \times_{\mathbb G_m} G)(k)$-linear isomorphism from $\omega_{V, W, \psi_k}$ to $\omega_{W, V, \psi_k}$. 
\end{itemize}
\subsubsection{The global Weil representation}
Now turn to the global situation, and take $k = \Q$ in (\ref{appendix bilinear forms subsubsec}).  
The adelic Schwartz space $\mathcal S_{\bA}(W_2\otimes V)$ is the restricted tensor product of the local Schwartz spaces $\mathcal S_{k}(W_2\otimes V)= \mathcal S(W_2\otimes V\otimes k)$ as $k$ ranges over completions of $\Q$. The global Weil representation $\omega = \omega_{V, W, \psi}$ of $(H \times_{\mathbb G_m} G)(\bA)$ is realized on $\mathcal S_{\bA}(W_2\otimes V)$ as the restricted tensor product of the local Weil representations.

Recall the automorphic realization of $\omega$, given by the theta kernel:
 \begin{equation}
     \theta (h,g;\phi) =\sum_{x\in W_2(\Q)\otimes V(\Q)}\omega(h,g)\phi (x),\;\; \;(h,g)\in (H\times_{\mathbb G_m} G)(\bA), \;\;\phi\in\mathcal S_{\bA}(W_2\otimes V).  
 \end{equation}

\subsubsection{Theta lifts of automorphic forms} \label{appendix global theta lift subsubsec}
Let $f\in\mathcal A_0 (G) $ be an automorphic cusp form and choose any $\phi\in\mathcal S_{\bA} (W_2\otimes  V)$. Then, fixing a Haar measure $\d g_1 $
 on $G_1 (\bA) $, the similitude theta lift $\theta_\phi(f)$  to $H$ is the automorphic function
 \begin{equation}\label{eq: glopbal theta lift def}
     h\mapsto \int_{[G_1]} \theta(g_1g_0,h; \phi) {f (g_1g_0)}\d h_1, \;\;\; h\in H({\bA}),
 \end{equation}
 where $g_0\in G({\bA})$ is any element such that $\nu_G(g_0) = \nu_H(h).$ 

 For any compact open subgroup $K \subset H(\bA_f)$, we say $\phi_f\in \mathcal S_{\bA_f} (W_2\otimes V)$ is $K$-invariant, which we write as $$\phi_f\in  \mathcal S_{\bA_f} (W_2\otimes V)^K,$$ if for all $k\in K$, there exists $g_0 \in G(\bA_f)$ with $\nu_G(g_0) = \nu_H(k)$ such that $$\omega(g_0, k) \phi_f = \phi_f.$$ Note that, if we fix $\phi_\infty\in \mathcal S_\R(W_2\otimes V)$, then 
 \begin{equation}
     \theta_{\phi_f\otimes\phi_\infty} (f)\textrm{ is $K$-invariant for all }\phi_f\in \mathcal S_{\bA_f} (W_2\otimes V)^K.
 \end{equation}
\subsection{Some explicit Schwartz functions}\label{appendix section explicit schwartz}
\subsubsection{The split four-dimensional quadratic space}
We briefly recall the conventions of \cite[\S5.1]{sweeting2022tateclassesendoscopyoperatornamegsp4}. 
Let $V = M_2$, with its canonical involution $x \mapsto x^\ast$ and quadratic form given by $(x, y) = \operatorname{tr}(xy^\ast).$ We have the map of algebraic groups over $\Q$:
\begin{equation}
    \boldsymbol{p}_Z: \GL_2 \times \GL_2 \to \operatorname{GO}(V)
\end{equation}
defined by $\boldsymbol{p}_Z(g_1, g_2) \cdot x = g_1 x g_2^\ast$. The kernel of $\boldsymbol{p}_Z$ is the antidiagonally embedded $\mathbb G_m$, and $\boldsymbol{p}_Z$ is a surjection onto the connected similitude group $\operatorname{GSO}(V) \subset \operatorname{GO}(V)$.

\subsubsection{}For any pair of automorphic forms $f_1,f_2$ on $\GL_2(\bA)$ with the same central character, we obtain an automorphic form $f_1\boxtimes f_2$ on $\operatorname{GSO}(V)(\bA)$ defined by \begin{equation}
    (f_1\boxtimes f_2)(\boldsymbol{p}_Z(g_1,g_2)) = f_1(g_1)f_2(g_2),\;\; g_1,g_2\in \GL_2(\bA).
\end{equation}
\subsubsection{Symplectic spaces}
For all $g$, we consider the standard symplectic space of dimension $2g$ over $\Q$, with basis $e_1,e_2,\ldots, e_{2g}$ such that
$$ e_{2n-1}\cdot e_{2n} = -e_{2n}\cdot e_{2n-1} = 1,\, \forall 1\leq n \leq g,$$ and all other pairings of basis vectors are trivial. We will always take the complete polarization $$\langle e_1, e_2, \ldots, e_{2n} \rangle =  \langle e_1,e_3,\ldots, e_{2g-1}\rangle\oplus  \langle e_2,e_4,\ldots, e_{2g}\rangle.$$
Note these are not the same coordinates as used in the main text; the change is to match with \cite{sweeting2022tateclassesendoscopyoperatornamegsp4}.

\subsubsection{Nonarchimedean Schwartz functions}
For each prime $p$, define the Schwartz function $$\phi_{N,p} \in \mathcal S_{\Q_p} (V)$$ to be the indicator function of the subset $$\begin{pmatrix}
    \Z_p & \Z_p \\ N\Z_p & \Z_p 
\end{pmatrix}\subset V\otimes \Q_p = M_2(\Q_p).$$
Clearly $\phi_{N,p}$ depends only on the $p$-adic valuation of $N$.
Also identify  $\mathcal S_{\Q_p} (V)$ with the Schwartz space $\mathcal S_{\Q_p} (\langle e_2\rangle\otimes V)$, which realizes the Weil representation of $(\GSp_2\times_{\mathbb G_m} \operatorname{GSO(V)})(\Q_p)$.

Fix the polarization $V = V_1\oplus V_2,$  where
 \begin{equation}\label{eq:isotropic of V}
   V_1 =   \begin{pmatrix}
        x & y\\0 & 0
      \end{pmatrix},\;\; V_2 =\begin{pmatrix}0 & 0\\z & w\end{pmatrix}.
  \end{equation}  
\begin{prop}\label{appendix prop:easy fourier transform}
    Under the Fourier transform of (\ref{appendix eqn:Fourier_transform}), $$\widehat \phi_{N,p} \in \mathcal S_{\Q_p} (\langle e_1,e_2\rangle\otimes V_2)$$
    is the indicator function of the set $$e_1 \otimes \begin{pmatrix} 0 & 0 \\ \Z_p & \Z_p\end{pmatrix}+ e_2\otimes \begin{pmatrix} 0 & 0 \\ N\Z_p & \Z_p\end{pmatrix}\subset \langle e_1,e_2\rangle\otimes V_2.$$
\end{prop}
\begin{proof}
    By definition, $$\widehat\phi_{N,p}\left(e_1 \otimes\begin{pmatrix}
        0 & 0 \\ z_1 & w_1
    \end{pmatrix} + e_2 \otimes\begin{pmatrix}
         0 & 0 \\ z_2 & w_2
    \end{pmatrix}\right) = \int_{\Q_p^2} \phi_{N,p}\left(\begin{pmatrix}
        x & y \\ z_2 & w_2 
    \end{pmatrix}\right) \cdot \psi \left( xw_1 - yz_1\right)\d x\d y, $$
    and the proposition follows. 
\end{proof}
\begin{prop}
    Fix integers $N$ and $M$ and consider the Schwartz function $$\phi_{N,M,p} \colonequals \phi_{N,p}\otimes \phi_{M,p} \in \mathcal S_{\Q_p} (\langle e_2 \rangle\otimes V)\otimes \mathcal S_{\Q_p} (\langle e_4 \rangle\otimes V)\subset \mathcal S_{\Q_p} (\langle e_2,e_4 \rangle\otimes V).$$ Then we have $$\phi_{N,M,p} \in  \mathcal S_{\Q_p} (\langle e_2,e_4 \rangle\otimes V)^{K_{N, 2,p}\cap K_{M,2,p}}.$$
    \end{prop}
\begin{proof}
    By the same calculation as Proposition \ref{appendix prop:easy fourier transform}, the Fourier transform $$\widehat\phi_{N,M,p} \in \mathcal S_{\Q_p} (\langle e_1,e_2,e_3,e_4\rangle \otimes V_2)$$
    is the indicator function of the set
    $$e_1 \otimes \begin{pmatrix}
        0 & 0 \\ \Z_p & \Z_p 
    \end{pmatrix} + e_2 \otimes \begin{pmatrix}
        0 & 0 \\ N\Z_p & \Z_p 
    \end{pmatrix} + e_3 \otimes \begin{pmatrix}
        0 & 0 \\ \Z_p & \Z_p 
    \end{pmatrix} + e_4 \otimes \begin{pmatrix}
        0 & 0 \\ M\Z_p & \Z_p 
    \end{pmatrix}.$$
    Because the Fourier transform is equivariant for $(\GSp_4\times_{\mathbb G_m}\operatorname{GSO}(V))(\Q_p)$, the proposition follows from the stability of this  set under the action of $$K_{N,2,p}\cap K_{M, 2,p} = \set{ g\in \GSp_4(\Z_p): g \equiv \operatorname{Id} \pmod{ p^{\operatorname{max}\set{\ord_p(N),\ord_p(M)}}}}.$$
\end{proof}
\subsubsection{The local Siegel-Weil map}
For each place $v$ of $\Q$, 
we have a map $$M_{1,v}[\cdot]: \mathcal S_{\Q_v}(V) \to I_v(1/2)$$ given by $$M_{1,v}[\phi](g) = \omega(h_{\det(g)}, \boldsymbol{p}_Z(g, 1)) \widehat \phi(0) =  |\det(g)|^{-1} \int_{\Q_v^2}\phi\left(g^{-1} \begin{pmatrix}
    x & y \\ 0 & 0 \end{pmatrix}\right)\d x\d y,
$$
cf. \cite[\S6.4.6]{sweeting2022tateclassesendoscopyoperatornamegsp4}.
\begin{prop}\label{appendix prop:identify siegel weil section}
Suppose $p |N$. Then $$M_{1,p}[\phi_{N,p} - p^{-1} \phi_{N/p,p}] = (1 - p^{-1}) \varphi^0_{N,p}.$$
\end{prop}
\begin{proof}
First, we calculate, for $i \geq 0$:
 \begin{align*}
        M_{1,p}[\phi_{N,p}]\left(\begin{pmatrix}
        1 & 0 \\ p^i & 1
    \end{pmatrix}\right) &= \int_{\Q_p^2} \phi_{N,p} \left( \begin{pmatrix}
        x & y \\ p^i x & p^i y 
    \end{pmatrix}\right) \d x\d y \\
    &= \operatorname{Vol}(\Z_p \cap p^{\ord_p(N)-i} \Z_p) = \begin{cases} p^{i-\ord_p(N)}, & i \leq \ord_p(N) \\ 1, & i > \ord_p(N).\end{cases}
    \end{align*}
    In particular, \begin{equation}\label{appendix eqn:compute local siegel weil section}
        M_{1,p}[\phi_{N,p} - p^{-1}\phi_{N/p,p}]\left(\begin{pmatrix}
        1 & 0 \\ p^i & 1 
    \end{pmatrix}\right)  = \begin{cases}
        1 - p^{-1}, &  i \geq \ord_p(N) \\ 0, & 0 \leq i < \ord_p(N).
    \end{cases}
    \end{equation}
On the other hand, it is clear that $M_{1,p}[\phi_{N,p} - p^{-1}\phi_{N/p,p}]$ is invariant under $K_1(N)_p$.  
    By the Iwasawa decomposition $\GL_2(\Q_p) = B(\Q_p) \cdot \GL_2(\Z_p)$ and the coset decomposition $$ \GL_2(\Z_p) = \bigsqcup_{0 \leq i \leq \ord_p(N)} \begin{pmatrix}
        1 & 0 \\ p^i & 1 
    \end{pmatrix} K_1(N)_p,$$ we conclude from (\ref{appendix eqn:compute local siegel weil section}) that $M_{1,p}[\phi_{N,p} - p^{-1}\phi_{N/p,p}]= (1- p^{-1})\varphi^0_{N,p}.$
\end{proof}

\subsubsection{Archimedean Schwartz function}
Let $\tau$ be the representation of $U(2)$ of highest weight $(3, -1)$.
We fix the nontrivial vector-valued archimedean Schwartz function
$$\phi_\infty \in \left( \mathcal S_{\R}(W_2\otimes V)\otimes \tau \otimes(\chi_2^\vee\boxtimes \chi_4^\vee)\right)^{U(2) \times \boldsymbol{p}_Z(\SO(2)\times \SO(2))},$$ denoted $\varphi^-_{4,2}$ in \cite[\S7.1.6]{sweeting2022tateclassesendoscopyoperatornamegsp4}. 

\subsection{Proof of Theorem \ref{main appendix theorem}}\label{appendix subsection proof}
\subsubsection{Construction of cohomology classes}
Fix new cuspidal Hecke eigenforms $f_1$ and $f_2$ for $\Gamma_1(N)$ of weights 4 and 2, respectively, and of equal nebetype character $\varepsilon$. Then  $f_1$ and $f_2$ correspond to automorphic forms
$$f_{1, \bA} \in \left(\mathcal A_0(\GSp_2; K_{N,1}) \otimes \chi_4\right)^{\R^\times \cdot U(1)},\;\; f_{2,\bA} \in \left(\mathcal A_0(\GSp_2; K_{N,1})\otimes \chi_2\right)^{\R^\times\cdot U(1)}.$$ 
For any Schwartz function $$\phi_f\in \mathcal S_{\bA_f}( \langle e_2, e_4 \rangle \otimes V)^{K_{N,2}},$$
we consider the vector-valued lift $$\Theta_{\phi_f\otimes \phi_\infty} (f_{1,\bA}\boxtimes f_{2,\bA}) \in \left(\mathcal A(\GSp_4; K_{N,2})\otimes \tau\right)^{\R^\times \cdot U(2)}.$$
\begin{rem}
    Assuming it is nonzero, the vector-valued automorphic form $\Theta_{\phi_f\otimes \phi_\infty} (f_{1,\bA}\boxtimes f_{2,\bA})$ generates the unique generic member of the endoscopic Yoshida lift $L$-packet on $\GSp_4$ associated to $f_1$ and $f_2$, cf. \cite{roberts2001global}. 
\end{rem}

By \cite[Theorem 8.3]{roberts2001global}, $\Theta_{\phi_f\otimes\phi_\infty} (f_{1,\bA}\boxtimes f_{2,\bA}) $ is a cusp form. Moreover, an easy calculation shows that $\operatorname{Hom}_{U(2)}(\tau, \wedge^3 \mathfrak p_\C^\ast)$ is one-dimensional in the notation of (\ref{appendix subsec:automorphic forms}) with $g = 2$. Hence from (\ref{appendix eqn:cusp form realization map}), we obtain a class
\begin{equation}
    \left[\Theta_{\phi_f\otimes \phi_\infty} (f_{1,\bA}\boxtimes f_{2,\bA})\right]\in H^3_c (\mathcal A_2'(N), \C)
\end{equation}
which is well-defined up to a scalar multiple.

When $g = 1$, the space $\mathfrak p^\ast_\C$ is a direct sum $\chi_2 \oplus \chi_2^\vee$ as a $U(1)$-module. In particular, $$\overline {f_{2,\bA}} \in (\mathcal A_0(\GSp_2; K_{N,1})\otimes \chi_2^\vee)^{\R^\times\cdot U(1)}$$ defines a class $$[\overline {f_{2,\bA}}] \in H^1_c(\mathcal A_1'(N), \C).$$
This is the usual cohomology class attached to the holomorphic modular form $f_2\otimes \varepsilon^{-1}$. 
By the K\"unneth formula, we also have the cohomology class
\begin{equation}
      [\overline {f_{2,\bA}}] \boxtimes \left[\Theta_{\phi_f\otimes \phi_\infty} (f_{1,\bA}\boxtimes f_{2,\bA})\right]\in H^4_c(\mathcal A_1'(N)\times \mathcal A_2'(N), \C).
\end{equation}
\begin{prop}\label{appendix prop:pairing as period}
    Up to a nonzero scalar depending on the normalizations, the Poincar\'e duality pairing
  $$\langle \left[\mathcal A_1'(N)\times_{\mu_N}\mathcal A_1'(N)\right],  [\overline {f_{2,\bA}}] \boxtimes \left[\Theta_{\phi_f\otimes {\phi_\infty}} (f_{1,\bA}\boxtimes f_{2,\bA})\right]\rangle\in H^8_c(\mathcal A_1'(N) \times \mathcal A_2'(N), \C)\simeq \C$$ is given by  $$\int_{[Z_H \backslash H]}\Theta_{\phi_f\otimes \phi_\infty} (f_{1,\bA}\boxtimes f_{2,\bA})(\iota(h_1,h_2)) \overline {f_{2,\bA}} (h_1) \d(h_1,h_2),$$
    where $H = \GSp_2 \times_{\mathbb G_m} \GSp_2$ is given the coordinates $(h_1,h_2)$ and $\iota: H \hookrightarrow \GSp_4$ is the standard embedding. 
\end{prop}
\begin{proof}
    See \cite[Proposition 7.2.4]{sweeting2022tateclassesendoscopyoperatornamegsp4}.
\end{proof}
\begin{lem}\label{appendix lemma:newform condition for nontriviality}
Suppose $N > 1$ is an integer such that there exist cuspidal newforms $f_1$ and $f_2$ for $\Gamma_1(N)$ of weights 4 and 2, respectively, of equal nebentype character $\varepsilon$. Then 
$$0 \neq  \left[\mathcal A_1'(N)\times_{\mu_N} \mathcal A_1'(N)\right] \in H^4(\mathcal A_1'(N) \times \mathcal A_2'(N),\Q).$$
\end{lem}
\begin{proof}
Without loss of generality, we may assume $f_1$ and $f_2$ are  Hecke eigenforms. Then
    Proposition \ref{appendix prop:pairing as period} reduces us to showing the nonvanishing of the period that appears therein, for some choice of $\phi_f\in \mathcal S_{\bA_f}(\langle e_2,e_4\rangle \otimes V)^{K_{N,2}}.$
    
We now fix the Schwartz function $\phi_f\in \mathcal S_{\bA_f} (\langle e_2,e_4\rangle\otimes V)$ to be of the form $\phi_{f}^{(1)}\otimes \phi_{f}^{(2)} $ for  $\phi_{f}^{(i)}= \otimes_p \phi_{p}^{(i)} \in \mathcal S_{\bA_f} (\langle e_{2i}\rangle \otimes V)$, $i = 1,2.$    By \cite[Theorem 6.5.2, Proposition 7.1.9]{sweeting2022tateclassesendoscopyoperatornamegsp4}, it suffices to show that, for all $p |N$, we may choose $\phi_{p}^{(i)}$ such that the following local zeta integrals 
are all nonzero:
\begin{equation}
\label{appendix eqn:local zeta integral}
\begin{split}
    \int_{(U\backslash \operatorname{PGL}_2\times U\backslash \operatorname{PGL}_2)(\Q_p)} \int_{\SL_2(\Q_p)} W_{\pi_{1,p}, \psi_{\Q_p}}^0(g_1) W_{\pi_{2,p},\psi_{\Q_p}}^0(g_2) W_{\pi_{2,p}^\vee,\psi_{\Q_p}^{-1}}^0(h_1h_c)\\ \omega(h_1h_c,g) \widehat \phi_{p}^{(1)}(1,0,0,-1) \varphi^0_{N,p}(g_2) M_{1,p}[\phi_{p}^{(2)}](g_1) \d h_1\d g_1\d g_2,\\
    c = \det (g_1g_2), \\g= \boldsymbol{p}_Z(g_1,g_2).
    \end{split}
\end{equation}
Here, $(1,0, 0, -1) \in \langle e_1,e_2\rangle \otimes V_2$ is the vector $e_1\otimes \begin{pmatrix} 0 & 0 \\ 1 & 0 \end{pmatrix}+ e_2 \otimes\begin{pmatrix}
     0 & 0 \\ 0 & -1
\end{pmatrix}$; $\pi_{1,p}$, $\pi_{2,p}$, and $\pi_{2,p}^\vee$ are the local components of the automorphic representations generated by $f_{1,\bA}$, $f_{2,\bA}$, and $\overline {f_{2,\bA}}$, respectively; and $W^0_{\pi_{1,p},\psi_{\Q_p}}$, etc. are the corresponding local newforms in the Whittaker models.  In fact, in \cite[Theorem 6.5.2]{sweeting2022tateclassesendoscopyoperatornamegsp4}, $\varphi^0_N$ is replaced with $\varphi^0 \colonequals \varphi^0_1$; however, the proof of \emph{loc. cit.} still applies, as long as Proposition \ref{appendix prop:nonzero residue} is used to replace the explicit calculation of the residue in \cite[Proposition 6.4.10]{sweeting2022tateclassesendoscopyoperatornamegsp4}. 
We choose our Schwartz functions as follows:
\begin{itemize}
    \item $\phi_p^{(1)} = \phi_{N,p}$ for all $p$.
    \item $\phi_p^{(2)} = \phi_{N,p} - p^{-1} \phi_{N/p,p} $ for all $p$, where $\phi_{N/p, p}$ is interpreted as 0 if $p \nmid N$.
\end{itemize}
With these choices, we now show that (\ref{appendix eqn:local zeta integral}) is nonzero.

We first consider the inner integral:
\begin{equation}\label{appendix eqn:inner integral}
   I(g) = \int_{\SL_2(\Q_p)} W_{\pi_{2,p}^\vee,\psi_{\Q_p}}^0(h_1h_c) \omega(h_1h_c,g) \widehat \phi_{N,p}(1,0,0,-1) \d h_1.
\end{equation}
Now, \cite[Lemma 6.3.3]{sweeting2022tateclassesendoscopyoperatornamegsp4} and its proof identifies (\ref{appendix eqn:inner integral}) with a function in the $\psi_{\Q_p}^{-1}\boxtimes \psi_{\Q_p}^{-1}$-Whittaker model of the representation $\pi_{2,p}^\vee\boxtimes \pi_{2,p}^\vee$ of $\GL_2(\Q_p)\boxtimes\GL_2(\Q_p)$. 
Because $\phi_{N,p}$ is clearly invariant by $\boldsymbol{p}_Z(K_1(N)_p\times K_1(N)_p)$, and because $\ord_p(N)$ is the conductor of $\pi_{2,p}$,  (\ref{appendix eqn:inner integral}) is a scalar multiple of the local newform $W_{\pi_{2,p}^\vee,\psi_{\Q_p}^{-1}}^0\boxtimes W_{\pi_{2,p}^\vee,\psi_{\Q_p}^{-1}}^0$.
To show this scalar multiple is nonzero, we evaluate 
\begin{equation}
    I(1) = \int_{\SL_2(\Q_p)} W_{\pi_{2,p}^\vee,\psi_{\Q_p}}^0(h_1) \omega(h_1, 1) \widehat \phi_{N,p} (1,0,0,-1) \d h_1.
\end{equation}
Now by Proposition \ref{appendix prop:easy fourier transform}, we can calculate directly that $\omega(h_1, 1) \widehat \phi_{N,p} (1,0,0,-1) $ is the indicator function of $K_1(N)_p$. Hence $$I(1) = \operatorname{Vol}(K_1(N)_p) \cdot W_{\pi_{2,p}^\vee,\psi_{\Q_p}}^0(1) \neq 0$$ by Proposition \ref{appendix prop:whittaker newform on torus}.
Hence, up to a nonzero scalar, (\ref{appendix eqn:local zeta integral}) becomes, after using Proposition \ref{appendix prop:identify siegel weil section}:
\begin{equation}
    \int_{(U\backslash \operatorname{PGL}_2\times U\backslash\operatorname{PGL}_2)(\Q_p)}W_{\pi_{2,p}^\vee, \psi_{\Q_p}^{-1}}^0(g_1) W_{\pi_{2,p}^\vee,\psi_{\Q_p}^{-1}}^0(g_2) W^0_{\pi_1,\psi_{\Q_p}} (g_1) W^0_{\pi_{2,p},\psi_{\Q_p}}(g_2) \varphi^0_{N,p} (g_1) \varphi^0_{N,p} (g_2) \d g_1 \d g_2.
\end{equation}
This factors into the product of the two integrals
\begin{equation}\label{appendix eqn:first factored local zeta integral}
    \int_{(U\backslash \operatorname{PGL}_2)(\Q_p)} W^0_{\pi_1,\psi_{\Q_p}} (g_1) W_{\pi_{2,p}^\vee, \psi_{\Q_p}^{-1}}^0(g_1)\varphi^0_{N,p} (g_1) \d g_1,
    \end{equation}\begin{equation}\label{appendix eqn:second factored local zeta integral}
        \int_{(U\backslash \operatorname{PGL}_2)(\Q_p)}W_{\pi_{2,p}^\vee,\psi_{\Q_p}^{-1}}^0(g_2) W^0_{\pi_2,p,\psi_{\Q_p}}(g_2)\varphi^0_{N,p} (g_2) \d g_2.
\end{equation}
Now, because $\varphi^0_{N,p}$ is supported on $B(\Q_p) K_1(N)_p$ by definition, and the local newforms are all $K_1(N)_p$-invariant, (\ref{appendix eqn:first factored local zeta integral}) and (\ref{appendix eqn:second factored local zeta integral}) become, up to nonzero scalars:
\begin{equation*}
\begin{split}
    \int_{\Q_p^\times} W^0_{\pi_1,\psi_{\Q_p}} \left(\begin{pmatrix}
         t & 0 \\ 0 & 1
    \end{pmatrix}\right) W_{\pi_{2,p}^\vee, \psi_{\Q_p}^{-1}}^0\left(\begin{pmatrix}
         t & 0 \\ 0 & 1
    \end{pmatrix}\right) |t| \d ^\times
t,\\   \int_{\Q_p^\times} W^0_{\pi_2,\psi_{\Q_p}} \left(\begin{pmatrix}
         t & 0 \\ 0 & 1
    \end{pmatrix}\right) W_{\pi_{2,p}^\vee, \psi_{\Q_p}^{-1}}^0\left(\begin{pmatrix}
         t & 0 \\ 0 & 1
    \end{pmatrix}\right) |t| \d ^\times t
    \end{split}
\end{equation*}
Now an easy computation using Proposition \ref{appendix prop:whittaker newform on torus} shows that both of these integrals are nonzero, which proves the lemma.
\end{proof}

\begin{thm}\label{main appendix theorem}
    For $N = 11$ and all $N \geq 13$, and for all $b: (\Z/N\Z)^2 \hookrightarrow (\Z/N\Z)^{2g}$, we have $$0 \neq [\widetilde{\mathsf{NL}}_{2,1}^b(N)]^+\in H^4(\mathcal A_1(N)\times \mathcal A_2(N), \Q).$$
\end{thm}
\begin{proof}

Combining Proposition \ref{appendix prop:rewrite NL class} with Lemmas \ref{appendix lemma:divisible ok} and  \ref{appendix lemma:newform condition for nontriviality}, it suffices to show that, for all such $N$, there  exists $N_0|N$ satisfying the following condition:
\begin{equation}\label{appendix eqn:condition exist newforms}
    \tag{$\ast$}\begin{split}\textrm{There exist cuspidal newforms $f_1$ and $f_2$ of weights 4 and 2, respectively,} \\ \textrm{for  $\Gamma_1(N_0)$, with equal central characters $\varepsilon$.}\end{split}
\end{equation}
First, we note that $p = 11$ and all primes $p > 13$ satisfy (\ref{appendix eqn:condition exist newforms}). Indeed, for such $p$ it is known that there exists a cuspidal newform $f$ of weight 2 for $\Gamma_0(p)$, cf. \cite[Proposition B.3]{halberstadt2002courbes}; then $f^2$ is a cuspidal modular form of weight 4 for $\Gamma_0(p)$, which is necessarily new because there are no cusp forms of weight 4 for $\SL_2(\Z)$. 

To exhibit more integers satisfying (\ref{appendix eqn:condition exist newforms}), we give the following table (sorted by prime factorization of $N_0$), in which all the data and labels are taken from \cite{lmfdb}.
\begin{center}
\begin{tabular}{|c|c|c|c|} 
 \hline 
 $N_0$ & $f_1$ & $f_2$ & $\varepsilon$\\
 \hline
 \hline
 $16 = 2^4$ & 16.4.e.a & 16.2.e.a & 16.e \\
$ 27 = 3^3$ & 27.4.a.a & 27.2.a.a & triv\\
 $25 = 5^2$ & 25.4.d.a & 25.2.d.a & 25.d\\
 $49 = 7^2$ & 49.4.a.a & 49.2.a.a & triv\\
  13  & 13.4.e.a & 13.2.e.a & 13.e\\
 $24 = 2^3\cdot 3$ & 24.4.a.a & 24.2.a.a & triv \\
 $18 = 2\cdot 3^2$ & 18.4.c.a & 18.2.c.a  & 18.c\\
$20 = 2^2\cdot 5$ & 20.4.a.a & 20.2.a.a & triv\\
$ 14 = 2 \cdot 7$& 14.4.a.a & 14.2.a.a & triv \\
 $15 = 3\cdot 5$ & 15.4.a.a & 15.2.a.a & triv\\
 $21= 3\cdot 7$ & 21.4.a.a & 21.2.a.a & triv \\
  $35 = 5\cdot 7$& 35.4.a.a & 35.2.a.a & triv\\
 \hline
\end{tabular}
\end{center}
Now a direct calculation shows that all $N$ as in the theorem are divisible by either 11, a prime $p > 13$, or one of the $N_0$ appearing in the table, which completes the proof. 
\end{proof}
    \bibliographystyle{plain} 
\bibliography{levelN_v6.bib}

\end{document}